%% file: prelie.tex
\documentclass[11pt,twoside,a4paper]{article}
\usepackage{makeidx}
\usepackage[dvips]{graphicx}
\usepackage{amscd}
\usepackage{amsmath}
\usepackage{amssymb}
\usepackage{latexsym}
\usepackage[latin1]{inputenc}
\usepackage[T1]{fontenc}   
\usepackage[english]{babel}

\input{xy}
\xyoption{all}

\newcommand{\M}{{\cal M}}

\newcommand{\DD}{{\cal DD}}
\newcommand{\g}{\mathfrak{g}}

\newcommand{\LI}{{\cal LI}}
\newcommand{\BI}{{\cal BI}}
\newcommand{\GM}{{\cal GM}}
\newcommand{\VGM}{{\cal VGM}}
\newcommand{\PL}{{\cal PL}}
\newcommand{\AS}{{\cal AS}}
\newcommand{\SI}{{\cal S}\Lambda}
\newcommand{\SIp}{{\cal S}_+\Lambda}
\newcommand{\tdelta}{\tilde{\Delta}}

\newcommand{\TI}{{\cal T}I}
\newcommand{\TIp}{{\cal T}_+I}
\newcommand{\U}{{\cal U}}

\title{Left ideals in an enveloping algebra, prelie products and applications to simple complex Lie algebras}
\author{L. Foissy \\
\\
{\small{\it Laboratoire de Mathématiques, Université de Reims}}\\
\small{{\it Moulin de la Housse - BP 1039 - 51687 REIMS Cedex 2, France}}\\
\small{e-mail : loic.foissy@univ-reims.fr}}
\date{}

\newtheorem{defi}{\indent Definition}
\newtheorem{lemma}[defi]{\indent Lemma}
\newtheorem{cor}[defi]{\indent Corollary}
\newtheorem{theo}[defi]{\indent Theorem}
\newtheorem{prop}[defi]{\indent Proposition}

\newenvironment{proof}{{\bf Proof.}}{\hfill $\Box$}

\begin{document}

\maketitle

\input{abstract.tex}

\tableofcontents

\input{intro.tex}
\input{chap1.tex}

\input{chap2.tex}

\input{chap3.tex}

\input{chap4.tex}

\bibliographystyle{amsplain}
\bibliography{biblio}

\end{document}

%% file: abstract.tex
ABSTRACT. We characterize prelie algebras in words of left ideals of the enveloping algebras
and in words of modules, and use this result to prove that a simple complex finite-dimensional
Lie algebra is not prelie, with the possible exception of $\mathfrak{f}_4$.\\



%% file: intro.tex
\section*{Introduction}

A (left) prelie algebra, or equivalently a left Vinberg algebra, or a left-symmetric algebra \cite{Chapoton1,Chapoton2,vanderLaan} is a couple $(V,\star)$,
where $V$ is a vector space and $\star$ a bilinear product on $V$ such that for all $x,y,z \in V$:
$$(x \star y)\star z-x \star (y \star z)=(y \star x)\star z-y \star (x \star z).$$
This axiom implies that the bracket defined by $[x,y]=x \star y-y\star x$ satisfies the Jacobi identity, so $V$ is a Lie algebra; 
moreover, $(V,\star)$ is a left module on $V$. The free prelie algebra is described in \cite{Chapoton1} in terms of rooted trees, 
making a link with the Connes-Kreimer Hopf algebra of Renormalization \cite{Connes,Kreimer}. \\

The aim of this text is to give examples of Lie algebras which are not prelie, namely the simple complex Lie algebras of finite dimension.
For this, we use the construction of \cite{Oudom} of the extension of the prelie product of a prelie algebra $(\g,\star)$ which gives to the symmetric
coalgebra  $S(\g)$ an associative product $*$, making it isomorphic to $\U(\g)$. This construction is also used in \cite{Loday2} 
to classify right-sided cocommutative graded and connected Hopf algebras.
A remarkable corollary of this construction is that $S(\g)_{\geq 2}=\bigoplus_{n\geq 2} S^n(\g)$ is a left ideal of $(S(\g),*)$. 
As a consequence, there exists a left ideal $I$ of $\U(g)$, such that the augmentation ideal $\U(\g)_+$ of $\U(g)$ can be decomposed as $\U(\g)_+=\g\oplus I$. 
We prove here the converse result: more precisely, given any Lie algebra $\g$, we prove in theorem \ref{11} that there exists a bijection into these two sets:
\begin{itemize}
\item $\PL(\g)=\{\star$ $\mid$ $\star$ is a  prelie product  on $\g$ inducing the bracket of $\g\}$.
\item $\LI(\g)=\{I$ left ideal of $\U(g)$ $\mid$ $\U_+(\g)=\g \oplus I\}$.
\end{itemize}

Let then a prelie algebra $\g$, and $I$ the ideal corresponding to the prelie product of $\g$. The $\g$-module $\U(\g)/I$ contains a submodule of codimension $1$,
isomorphic to $(\g,\star)$, and the special element $1+I$. This leads to the definition of good pointed modules (definition \ref{13}). 
The sets $\PL(\g)$ and $\LI(\g)$ are in bijection with the set $\GM(\g)$ of isoclasses of good pointed modules, as proved in theorem \ref{14}.

If $\g$ is semisimple, then any good pointed module is isomorphic to $K \oplus (\g,\star)$ as a $\g$-module. This leads to the definition
of very good pointed module, and we prove that if $\g$ is a semisimple Lie algebra, then it is prelie if, and only if $\g$ has a very good pointed module.
If $\g$ is a simple complex lie algebra, with the possible exception of $\mathfrak{f}_4$, we prove that it has no very good pointed module, so $\g$ is not prelie.
The proof is separated into three cases: first, the generic cases, then $\mathfrak{so}_{2n+1}$, and finally $\mathfrak{sl}_{6}$ and $\mathfrak{g}_{2}$, which 
are proved by direct computations using the computer algebra system MuPAD pro 4. We conjecture that this is also true for $\mathfrak{f}_{4}$ 
(the computations would be similar, though very longer). \\

We also give in this text a non cocommutative version of theorem \ref{11}, replacing enveloping algebras (which are symmetric coalgebras, as the base field
is of characteristic zero) by cofree coalgebras, and prelie algebras by dendriform algebras \cite{Loday1,Loday3}. 
If $A$ is a cofree coalgebra, we prove in theorem \ref{32} that there is a bijection between these two sets:
\begin{itemize}
\item $\DD(A)=\{(\prec,\succ)$ $\mid$ $(A,\prec,\succ,\Delta)$ is a dendriform Hopf algebra$\}$.
\item $\LI(A)=\left\{(*,I)\mid \begin{array}{l}
(A,*,\Delta)\mbox{ is a Hopf algebra and}\\
I\mbox{ is a left ideal of }A \mbox{ such that }A_+=Prim(A) \oplus I
\end{array}\right\}$.
\end{itemize}

The text is organized as follows: the first section deals with the general results on prelie algebras: after preliminaries on symmetric coalgebras,
theorem \ref{11} is proved and good pointed modules are introduced. The second section is devoted to the study of good pointed modules
over simple complex Lie algebras. The results on dendriform algebras are exposed in the third section and the MuPAD procedures used 
in this text are written in the last section. \\

{\bf Notations.} \begin{enumerate}
\item $K$ is a commutative field of characteristic zero. Any algebra, coalgebra, bialgebra, etc, of this text will be taken over $K$.
\item Let $\g$ be a Lie algebra. We denote by $\U(\g)$ its enveloping algebra and by $\U_+(\g)$ the augmentation ideal of $\U(\g)$.
\end{enumerate}

%% file: chap1.tex
\section{Prelie products on a Lie algebra}

\subsection{Preliminaries and recalls on symmetric coalgebras}

Let $V$ be a vector space. The algebra $S(V)$ is given a coproduct $\Delta$, defined as the unique algebra morphism from $S(V)$ to 
$S(V) \otimes S(V)$, such that $\Delta(v)=v\otimes 1+1\otimes v$ for all $v \in V$. Let us recall the following facts:
\begin{itemize}
\item Let us fix a basis $(v_i)_{i\in \Lambda}$ of $V$. We define $\SI$ as the set of sequences $a=(a_i)_{i\in \Lambda}$ of elements of $\mathbb{N}$,
with a finite support. For an element $a \in S\Lambda$, we put $\displaystyle l(a)=\sum_{i\in \Lambda} a_i$. For all $a\in \SI$, we put:
$$v_a=\prod_{i\in \Lambda} \frac{v_i^{a_i}}{a_i!}.$$
Then $(v_a)_{a\in \SI}$ is a basis of $S(V)$, and the coproduct is given by:
$$\Delta(v_a)=\sum_{b+c=a}v_b \otimes v_c.$$
\item Let $S_+(V)$ be the augmentation ideal of $S(V)$. It is given a coassociative, non counitary coproduct $\tdelta$ defined by
$\tdelta(x)=\Delta(x)-x\otimes 1-1\otimes x$ for all $x \in S_+(V)$. In other terms, putting $\SIp=\SI-\{(0)\}$, 
$(v_a)_{a\in \SIp}$ is a basis of $S_+(V)$ and:
$$\tdelta(v_a)=\sum_{\substack{b+c=a\\b,c \in \SIp}} v_b \otimes v_c.$$
\item Let $\tdelta^{(n)}:S_+(V) \longrightarrow S_+(V)^{\otimes (n+1)}$ be the $n$-th iterated coproduct of $S_+(V)$. Then:
$$Ker\left(\tdelta^{(n)}\right)=\bigoplus_{k=1}^n S^k(V).$$
In particular, $Prim(S(V))=V$.
\item $(S^n(V))_{n\in \mathbb{N}}$ is a gradation of the coalgebra $S(V)$.
\end{itemize}

\begin{lemma} \label{1}
In $S_+(V)\otimes S_+(V)$, $Ker(\tdelta \otimes Id-Id \otimes \tdelta)=Im(\tdelta)+V \otimes V$.
\end{lemma}

\begin{proof} $\supseteq$. Indeed, if $y=\tdelta(x)+v_1 \otimes v_2 \in Im(\tdelta)+V \otimes V$, as $\tdelta$ is coassociative:
\begin{eqnarray*}
(\tdelta \otimes Id)(y)-(Id \otimes \tdelta)(y)&=&(\tdelta \otimes Id)\circ \tdelta(x)-(Id \otimes \tdelta) \circ \tdelta(x)\\
&&+\tdelta(v_1)\otimes v_2-v_1 \otimes \tdelta(v_2)\\
&=&0.
\end{eqnarray*}

$\subseteq$. Let $X \in Ker(\tdelta \otimes Id-Id \otimes \tdelta)$. Choosing a basis $(v_i)_{i\in \Lambda}$ of $V$, we put:
$$X=\sum_{a,b \in \SIp}x_{a,b} v_a \otimes v_b.$$
By homogeneity, we can suppose that $X$ is homogeneous of a certain degree $n\geq 2$ in $S_+(V) \otimes S_+(V)$.
If $n=2$, then $X \in V \otimes V$. Let us assume that $n \geq 3$. Then:
$$(\tdelta \otimes Id) \circ \tdelta(X)=\sum_{a,b,c \in \SIp} x_{a+b,c} v_a \otimes v_b \otimes v_c
=(Id \otimes \tdelta) \circ \tdelta(X)=\sum_{a,b,c \in \SIp} x_{a,b+c} v_a \otimes v_b \otimes v_c.$$
So, for all $a,c \in \SIp$, $b\in \SI$, $x_{a+b,c}=x_{a,b+c}$.\\

Let $a,b,a',b' \in \SIp$, such that $a+b=a'+b'$. Let us show that $x_{a,b}=x_{a',b'}$. 

{\it First case.} Let us assume that the support of $a$ and $a'$ are not disjoint. For all $i \in \Lambda$, we put:
\begin{eqnarray*}
c_i&=&\left\{ \begin{array}{rcl}
a'_i-a_i&\mbox{if}& a'_i-a_i>0,\\
0&\mbox{if}&a'_i-a_i\leq 0,
\end{array}\right.\\
c'_i&=&\left\{ \begin{array}{rcl}
a_i-a'_i&\mbox{if}& a_i-a'_i>0,\\
0&\mbox{if}&a_i-a'_i\leq 0.
\end{array}\right. \end{eqnarray*}
Then $c$ and $c'$ belong to $\SI$. Moreover, for all $i\in \Lambda$, $a_i+c_i=a'_i+c'_i$, so $a+c=a'+c'$, or $a'=a+c-c'$. As $a+b=a'+b'$, $b'=b-c+c'$.

For all $i\in \Lambda$:
$$a_i-c'_i= \left\{ \begin{array}{rcl}
a'_i&\mbox{if}& a_i-a'_i>0,\\
a_i&\mbox{if}&a_i-a'_i\leq 0.
\end{array}\right.$$
So $a-c'\in \SI$. Moreover, if $a-c'=0$, then if $a_i>a'_i$, $a'_i=0$, so the support of $a'$ is included in $\{i\:/\:a_i \leq a'_i\}$.
If $a_i \leq a'_i$, then $a_i=0$, so the support of $a$ is included in $\{i\:/\: a_i>a'_i\}$. As a consequence, the supports of $a$ and $a'$
are disjoint: this is a contradiction. So $a-c' \in \SIp$. As $a-c'$ and $b' \in \SIp$: 
$$x_{a',b'}=x_{a-c'+c,b'}=x_{a-c',b'+c}=x_{a-c',b+c'}.$$
As $a-c'$ and $b \in \SIp$:
$$x_{a-c',b+c'}=x_{a,b}.$$
The proof is similar if the support of $b$ and $b'$ are not disjoint, permuting the roles of $(a,a')$ and $(b,b')$.\\

{\it Second case.} Let us assume that the supports of $a$ and $a'$, and the supports of $b$ and $b'$ are disjoint. 
We then denote, up to a permutation of the index set $\Lambda$:
\begin{eqnarray*}
a&=&(a_1,\ldots,a_k,0,\ldots,0,\ldots),\\
a'&=&(0,\ldots,0,a'_{k+1},\ldots,a'_{k+l},0,\ldots),
\end{eqnarray*}
where the $a_i$'s and the $a'_j$'s are non-zero. As $a+b=a'+b'$, necessarily $b_{k+1}\neq 0$.\\

{\it First subcase.} $l(a)$ or $l(a')>1$. We can suppose that $l(a)>1$. Then $a$ and $(a_1-1,a_2,\ldots, a_k,0,\ldots)$ 
are elements of $\SIp$ because $l(a)>1$ and their support are non disjoint. By the first case:
$$x_{a,b}=x_{(a_1-1,a_2,\ldots,a_k,0,\ldots),b+(1,0,\ldots)}.$$
Moreover, the supports of  $(a_1-1,a_2,\ldots, a_k,0,\ldots)$  and $(a_1-1,a_2,\ldots, a_k,1,0,\ldots)$ are not disjoint.
As $b_{k+1}>0$, $b+(1,0,\ldots)-(0,\ldots,0,1,0,\ldots)$ belongs to $\SIp$. By the first case:
$$x_{(a_1-1,a_2,\ldots,a_k,0,\ldots),b+(1,0,\ldots)}=x_{(a_1-1,a_2,\ldots,a_k,1,0,\ldots),b+(1,0,\ldots,0,-1,0,\ldots)}.$$
Finally, the support of $(a_1-1,a_2,\ldots,a_k,1,0,\ldots)$ and $a'$ are not disjoint, so:
$$x_{(a_1-1,a_2,\ldots,a_k,1,0,\ldots),b+(1,0,\ldots,0,-1,0,\ldots)}=x_{a',b'}.$$
If  $l(b)$ or $l(b')>1$, the proof is similar, permuting the roles of $(a,a')$ and $(b,b')$.\\

{\it Second subcase.} $l(a)=l(a')=l(b)=l(b')=1$. Then $v_a\otimes v_b, v_{a'} \otimes v_{b'} \in V \otimes V$.
By the homogeneity condition, $x_{a,b}=x_{a',b'}=0$.\\
 
Hence, we put $x_{a+b}=x_{a,b}$ for all $a,b \in \SIp$: this does not depend of the choice of $a$ and $b$. So:
$$X=\sum_{c \in \SIp} x_c \sum_{\substack {a+b=c\\a,b \in \SIp}}v_a \otimes v_b=\sum_{c \in \SIp} x_c \tdelta(v_c),$$
so $X \in Im(\tdelta)$. \end{proof}

\begin{lemma} \label{2}
In $S(V)$, $Im(\tdelta)\cap (V \otimes V)=S^2(V)$.
\end{lemma}

\begin{proof} $\subseteq$. By cocommutativity of $\tdelta$. $\supseteq$. If $v,w \in V$, then $v \otimes w+w\otimes v=\tdelta(vw)$. \end{proof}

\begin{lemma} \label{3}
Let $W$ be a subspace of $S_+(V)$, such that $S_+(V)=V \oplus W$. There exists a unique coalgebra endomorphism $\phi$ of $S(V)$,
such that $\phi_{\mid V}=Id_V$ and $\displaystyle \phi\left(\bigoplus_{n\geq 2} S^n(V)\right)=W$. Moreover, $\phi$ is an automorphism.
\end{lemma}

We shall denote from now:
$$S_{\geq 2}(V)=\bigoplus_{n=2}^\infty S^n(V).$$

\begin{proof} {\it Existence.} We denote by $\pi_W$ the projection on $W$ in the direct sum $S(V)=(1) \oplus V \oplus W$.
We define inductively $\phi_{\mid S^n(V)}$. If $n=0$, it is defined by $\phi(1)=1$. If $n=1$, it is defined by $\phi_{\mid V}=Id_V$. 
Let is assume that $\phi$ is defined on the subcoalgebra $C_{n-1}=K\oplus V \oplus \ldots \oplus S^{n-1}(V)$ with $n \geq 2$
and let us define $\phi$ on $S^n(V)$. Let $v_a \in S^n(V)$. Then $\tdelta(v_a) \in C_{n-1} \otimes C_{n-1}$,
so $(\phi \otimes \phi)\circ \tdelta(v_a)$ is already defined. Moreover:
\begin{eqnarray*}
(\tdelta\otimes Id)\circ (\phi \otimes \phi)\circ \tdelta(v_a)&=& (\phi \otimes \phi\otimes \phi)\circ(\tdelta\otimes Id)\circ \tdelta(v_a)\\
&=&(\phi \otimes \phi\otimes \phi)\circ(Id \otimes \tdelta)\circ \tdelta(v_a)\\
&=&(Id\otimes \tdelta)\circ (\phi \otimes \phi)\circ \tdelta(v_a).
\end{eqnarray*}
By lemma \ref{1}, $(\phi \otimes \phi)\circ \tdelta(v_a)\in Im(\tdelta) +V \otimes V$. Moreover, as $\tdelta$ is cocommutative, using lemma \ref{2}:
\begin{eqnarray*}
(\phi \otimes \phi)\circ \tdelta(v_a)&\in&(Im(\tdelta) +V \otimes V)\cap S^2(S(V))\\
&\in&Im(\tdelta)+(V \otimes V)\cap S^2(S(V))\\
&\in&Im(\tdelta)+S^2(V)\\
&\in&Im(\tdelta).
\end{eqnarray*}
Let $w_a \in S_+(V)$, such that $(\phi \otimes \phi)\circ \tdelta(v_a)=\tdelta(w_a)$. We put then $\phi(v_a)=\pi_W(w_a)$.
So $\phi(v_a) \in W$. Moreover, $w_a-\pi_W(w_a)\in V\subseteq Ker(\tdelta)$, so:
\begin{eqnarray*}
(\phi \otimes \phi) \circ \Delta(v_a)&=&\phi(v_a)\otimes \phi(1)+\phi(1)\otimes \phi(v_a)+(\phi \otimes \phi)\circ \tdelta(v_a)\\
&=&\phi(v_a)\otimes 1+1\otimes \phi(v_a)+\tdelta(w_a)\\
&=&\phi(v_a)\otimes 1+1\otimes \phi(v_a)+\tdelta(\pi_W(w_a))\\
&=&\Delta(\phi(v_a)).
\end{eqnarray*}
So $\phi_{\mid C_n}$ is a coalgebra morphism.\\

{\it Unicity.} Let $\tilde{\phi}$ be another coalgebra endomorphism satisfying the required properties. Let us show that $\phi(v_a)=\tilde{\phi}(v_a)$ 
by induction on $n=l(a)$. If $n=0$ or $1$, this is immediate. Let us assume the result at all rank $<n$, $n \geq 2$. Then:
$$\tdelta(\phi(v_a)-\tilde{\phi}(v_a))=(\phi \otimes \phi-\tilde{\phi} \otimes \tilde{\phi})\left(\sum_{b+c=a \in \SIp} v_b \otimes v_c \right)=0,$$
by the induction hypothesis. So $\phi(v_a)-\tilde{\phi}(v_a)\in Prim(S(V))=V$. Moreover, it belongs to 
$\phi(S_{\geq 2}(V))+\tilde{\phi}(S_{\geq 2}(V))\subseteq W$. As $V\cap W=(0)$, $\phi(v_1\ldots v_n)=\tilde{\phi}(v_1\ldots v_n)$. \\

We have defined in this way an endomorphism $\phi$, such that $\phi_{\mid V}=Id_V$ and $\phi(S^n(V))\subseteq W$ for all $n\geq 2$. 
We now show that $\phi$ is an automorphism. Let us suppose that $Ker(\phi)\neq (0)$. As $\phi$ is a coalgebra morphism, 
$Ker(\phi)$ is a non-zero coideal, so it contains primitive elements, that is to say elements of $V$: impossible, as $\phi_{\mid V}$ is monic. 
So $\phi$ is monic. Let us prove that $v_a\in Im(\phi)$ by induction on $l(a)$. If $l(a)=0$ or $1$, then $\phi(v_a)=v_a$. We can suppose that
$v_a=\lambda v_1\ldots v_k$, where $\lambda$ is a non-zero scalar. Then: 
$$\tdelta^{k-1}(v_a)=\lambda \sum_{\sigma\in S_k} v_{\sigma(1)} \otimes \ldots \otimes v_{\sigma(k)}
=\lambda \sum_{\sigma\in S_k} \phi(v_{\sigma(1)}) \otimes \ldots \otimes \phi(v_{\sigma(k)})=\tdelta^{k-1}(\phi(v_a)).$$
So $v_a-\phi(v_a) \in Ker(\tdelta^{k-1})=K\oplus V \oplus \ldots \oplus S^{k-1}(V)\subseteq Im(\phi)$ by the induction hypothesis.
So $v_a \in Im(\phi)$ and $\phi$ is epic. As a consequence, $\phi(S_{\geq 2}(V))=W$. \end{proof}

\subsection{Extension of a prelie product}

\begin{defi}\textnormal{
A (left) prelie algebra is a vector space $\g$, with a product $\star$ satisfying the following property: for all $x,y,z \in \g$,
$$(x \star y)\star z-x \star (y \star z)=(y \star x)\star z-y \star (x \star z).$$
}\end{defi}

{\bf Remark.} A prelie algebra $\g$ is also a Lie algebra, with the bracket given by:
$$[x,y]=x \star y-y\star x.$$ 
This bracket will be called the Lie bracket induced by the prelie product.\\

Let $(\g,\star) $ be a prelie algebra. By \cite{Gan} and \cite{Oudom}, permuting left and right, the prelie product $\star$ can be extended to the coalgebra $S(\g)$ 
in the following way: for all $x,y\in \g$, $P,Q,R \in S(\g)$,
$$\left\{ \begin{array}{rcl}
1\star P&=&P,\\
(xP)\star y&=&x \star (P \star y)-(x \star P) \star y,\\
P \star (QR)&=&\sum \left(P^{(1)} \star Q\right)\left(P^{(2)} \star R\right).
\end{array}\right.$$
Then $S(g)$ is given an associative product $*$ defined by $P*Q=\sum P^{(1)} \left(P^{(2)} \star Q\right)$. Moreover, $(S(\g),*,\Delta)$ is isomorphic, 
as a Hopf algebra, to $\U(\g)$. There exists a unique isomorphism $\xi_\star$ of Hopf algebras:
$$\xi_\star: \left\{ \begin{array}{rcl}
\U(\g)&\longrightarrow & (S(\g),*,\Delta)\\
x \in \g&\longrightarrow&x \in \g.
\end{array} \right.$$

\begin{lemma} \label{5}
For all $n \geq 1$, $\g \star S^n(\g) \subseteq S^n(\g)$.
\end{lemma}

\begin{proof} By induction on $n$. This is immediate for $n=1$. Let us assume the result at rank $n-1$. 
Let $P=yQ \in S^n(\g)$, with $y\in \g$ and $Q\in S^{n-1}(\g)$. For all $x\in \g$, as $x$ is primitive:
$$x \star P=(x \star y)Q+y (x \star Q).$$
Note that $x\star y \in \g$ and $x \star Q \in S^{n-1}(\g)$ by the induction hypothesis. So $x \star P \in S^n(\g)$. \end{proof}

\begin{prop} \label{6}
Let $\g$ be a prelie algebra. We denote $\displaystyle S_{\geq 2}(\g)=\bigoplus_{n\geq 2}S^n(\g)$. Then:
\begin{enumerate}
\item $S_{\geq 2}(\g)$ is a left ideal for $*$.
\item $S_{\geq 2}(\g)$ is a bilateral ideal for $*$ if, and only if, $\star$ is associative on $\g$.
\end{enumerate}
\end{prop}

\begin{proof} 1. Let $x\in \g$ and $Q\in S^n(\g)$, with $n\geq 2$. Then, by lemma \ref{5}:
$$x*Q=xQ+x \star Q \in S^{n+1}(V)+S^n(V).$$
As a consequence, $S_{\geq 2}(\g)$ is stable by left multiplication by an element of $\g$. As $\g$ generates $(S(\g),*)$
(because it is isomorphic to $\U(\g)$), $S_{\geq 2}(\g)$ is a left ideal.\\

2, $\Longleftarrow$. Let us first show that $S^n(\g) \star \g \subseteq S_{\geq 2}(\g)$ for all $n \geq 2$ by induction on $n$.
For $n=2$, take $x,y,z \in \g$. Then $(xy) \star z=x \star (y\star z)-(x\star y)\star z=0$,
as $\star$ is associative on $\g$. Let us assume the result at rank $n-1$ $(n \geq 3$). Let $x\in \g$, $P\in S^{n-1}(V)$, $y\in \g$. Then
$(xP) \star z=x \star (P\star z)-(x\star P)\star z$. By the induction hypothesis, $P\star z \in S_{\geq 2}(\g)$. 
By lemma \ref{5}, $x\star (P\star z) \in S_{\geq 2}(\g)$. By lemma \ref{5}, $x\star P \in S^{n-1}(\g)$. By the induction hypothesis, 
$(x \star P)\star z \in S_{\geq 2}(\g)$. So $S^n(\g) \star \g \subseteq S_{\geq 2}(\g)$ for all $n\geq 2$.

Let us now prove that $S_{\geq 2}(\g)*\g \subseteq S_{\geq 2}(\g)$. Let $P \in S_{\geq 2}(\g)$ and $y\in \g$.
We put $\Delta(P)=P\otimes 1+1\otimes P+\sum P'\otimes P''$. Then:
$$P*y=Py+P\star y+\sum P'(P''\star y).$$
Note that $Py$ and $\sum P'(p''\star y)$ belong to $S_{\geq 2}(\g)$. We already proved $P\star y \in S_{\geq 2}(\g)$.
As $\g$ generates $(S(\g),*)$, $S_{\geq 2}(\g)$ is a right ideal. By the first point, it is a bilateral ideal.\\

2, $\Longrightarrow$. Let us assume that $\star$ is not associative on $\g$. There exists $x,y,z\in \g$, 
such that $x \star (y\star z)-(x\star y)\star z \neq 0$. Then:
\begin{eqnarray*}
(xy)*z&=&xyz+x (y\star z)+y(x\star z)+(xy) \star z\\
&=&\underbrace{xyz+x (y\star z)+y(x\star z)}_{\in S_{\geq 2}(\g)}+\underbrace{x \star (y\star z)-(x\star y)\star z}_{\in V-\{0\}},
\end{eqnarray*}
so $(xy)* z \notin S_{\geq 2}(\g)$, which is not a right ideal. \end{proof}

\begin{defi}
\textnormal{Let $\g$ be a Lie algebra. We define:
\begin{enumerate}
\item $\PL(\g)=\{\star$ $\mid$ $\star$ is a  prelie product  on $\g$ inducing the bracket of $\g\}$.
\item $\LI(\g)=\{I$ left ideal of $\U(g)$ $\mid$ $\U_+(\g)=\g \oplus I\}$.
\end{enumerate}}
\end{defi}

\begin{prop}
There exists an application:
$$\Phi_\g: \left\{ \begin{array}{rcl}
\PL(\g) &\longrightarrow& \LI(\g)\\
\star &\longrightarrow & \xi_\star^{-1}(S_{\geq 2}(\g))\\
\end{array} \right.$$
\end{prop}

\begin{proof}  $\Phi_\g(\star)$ is indeed an element of $\LI(\g)$, as $S_{\geq 2}(\g)$ is a left ideal of $(S(\g),*)$
and $\xi_\star$ is an isomorphism of algebras. \end{proof}

\subsection{Prelie product associated to a left ideal}

\begin{prop}
Let $I \in \LI(\g)$. We denote by $\varpi_\g$ the projection on $\g$ in the direct sum $\U(\g)=(1)\oplus \g \oplus I$. Then the product $\star$
defined by $x \star y=\varpi_\g(xy)$ is an element of $\PL(\g)$. This defines an application $\Psi_\g:\LI(\g) \longrightarrow \PL(\g)$.
\end{prop}

\begin{proof} Let $x,y \in \g$. In $\U(\g)$, $xy-yx=[x,y]\in \g$, so:
$$x\star y-y\star x=\varpi_\g(xy-yx)=\varpi_\g([x,y])=[x,y],$$
so $\star$ induces the Lie bracket of $\g$. It remains to prove that it is prelie. Let us fix a basis $(v_i)_{i\in \Lambda}$ of $\g$. By the Poincaré-Birkhoff-Witt
theorem, $\U(\g)$ is isomorphic to $S(\g)$ as a coalgebra. Using lemma \ref{3},  there exists an isomorphism of coalgebras:
$$\phi_I: \left\{ \begin{array}{rcl}
S(\g)&\longrightarrow & \U(\g)\\
v \in \g&\longrightarrow &\g,
\end{array}\right.$$
such that $\phi_I(S_{\geq 2}(\g))=I$. We denote by $(v^a)_{a\in \SI}$ the basis of $\U(\g)$, image of the basis $(v_a)_{a\in \SI}$ of $S(\g)$. 
As a consequence:
\begin{itemize}
\item A basis of $I$ is given by $(v^a)_{a\in \SI, \: l(a)\geq 2}$.
\item For all $a\in \SI$, $\displaystyle \Delta(v^a)=\sum_{b+c=a}v^b \otimes v^c$.
\end{itemize}
We put $\delta_j=(\delta_{i,j})_{i\in I}$ for all $j\in I$ (so $v_i=v^{\delta_i}$ for all $i\in I$) and, for all $a,b \in \SI$:
$$v^av^b=\sum_{c\in \SI} x_{a,b}^c v^c.$$
By definition, for all $i,j \in I$:
$$v_i \star v_j=\sum_{k\in I} x_{\delta_i,\delta_j}^{\delta_k} v_k.$$

Let us prove the prelie relation for $v_i,v_j,v_k$. It is obvious if $i=j$: let us assume that $i\neq j$. Then, in $\U(\g)$:
\begin{eqnarray*}
\Delta\left(v_iv_j-v^{\delta_i+\delta_j}\right)&=&\left(v_iv_j-v^{\delta_i+\delta_j}\right)\otimes 1+1\otimes \left(v_iv_j-v^{\delta_i+\delta_j}\right)\\
&&+v_i \otimes v_j+v_j \otimes v_i-v^{\delta_i} \otimes v^{\delta_j}-v^{\delta_j} \otimes v^{\delta_i}\\
&=&\left(v_iv_j-v^{\delta_i+\delta_j}\right)\otimes 1+1\otimes \left(v_iv_j-v^{\delta_i+\delta_j}\right).
\end{eqnarray*}
So $v_iv_j-v^{\delta_i+\delta_j} \in Prim(S(\g))=\g$. So:
$$v_iv_j=v^{\delta_i+\delta_j}+\sum_{k\in I} v_{\delta_i,\delta_j}^{\delta_k} v^{\delta_k}=v^{\delta_i+\delta_j}+v_i\star v_j.$$
Hence:
\begin{eqnarray*}
\varpi_\g(v_i(v_jv_k))&=&\varpi_\g\left(\sum_{c\in \SI} c_{\delta_j,\delta_k}^c v_iv^c\right)\\
&=&\varpi_\g\left(\sum_{\substack{c\in \SI\\l(c)\geq 2}} c_{\delta_j,\delta_k}^c \underbrace{v_iv^c}_{\in I,\mbox{\scriptsize{ left ideal}}}\right)
+\varpi_\g(v_i(v_j \star v_k))\\
&=&0+v_i \star (v_j \star v_k)\\
&=&\varpi_\g((v_iv_j)v_k)\\
&=&\varpi_\g(v^{\delta_i+\delta_j}v_k)+\varpi_\g((v_i \star v_j)v_k)\\
&=&\varpi_\g(v^{\delta_i+\delta_j}v_k)+(v_i \star v_j)\star v_k.
\end{eqnarray*}
So:
$$v_i \star (v_j \star v_k)-(v_i \star v_j)\star v_k=\varpi_\g(v^{\delta_i+\delta_j}v_k)=v_j \star (v_i \star v_k)-(v_j \star v_i)\star v_k.$$
Hence, the product $\star$ is prelie. \end{proof}

\subsection{Prelie products on a Lie algebra}

\begin{prop}
The applications $\Phi_\g$ and $\Psi_\g$ are inverse bijections.
\end{prop}

\begin{proof} Let $\star \in \PL(\g)$. We put $I=\Phi_\g(\star)$ and $\bullet=\Psi_\g(I)$. Let $\pi_\g$ be the canonical surjection on $\g$ in $S(\g)$. 
Then, as $\xi_\star(I)=S_{\geq 2}(\g)$, $\xi_\star \circ \varpi_\g=\pi_\g \circ \xi_\star$. So, for all $x,y \in \g$:
$$\xi_\star(x \bullet y)=\xi_\star\circ \varpi_\g(xy)=\pi_\g\circ \xi_\star(xy)=\pi_\g(x*y)=\pi_\g(xy+x\star y)=x\star y=\xi_\star(x\star y).$$
As $\xi_\star$ is monic, $x\bullet y=x\star x$, so $\Psi_\g\circ \Phi_\g(\star)=\star$.

Let $I\in \LI(\g)$. We put $\star=\Psi_\g(I)$. We have to prove that $\xi_\star^{-1}(S_{\geq 2}(\g))=I$, that is to say $\xi_\star(I)=S_{\geq 2}(\g)$.
Because they are both complements of $(1)\oplus \g$ and $\xi_\star$ is bijective, it is enough to prove $\xi_\star(I)\subseteq S_{\geq 2}(\g)$.
With the preceding notations, let $v^a\in I$, $l(a)\geq 2$, and let us prove that $\xi_\star(v^a) \in S_{\geq 2}(\g)$. If $l(a)=2$, 
we put $a=\delta_i+\delta_j$. Then we saw that $v^a=v_iv_j-v_i \star v_j$ if $i\neq j$. If $i=j$, in the same way, $2v^a=v_iv_j-v_i \star v_j$.
So, up to a non-zero mutiplicative constant $\lambda$:
$$\xi_\star(v^a)=\lambda(v_i*v_j-v_i \star v_j)=\lambda(v_iv_j+v_i \star v_j-v_i \star v_j)=\lambda v_iv_j \in S_{\geq 2}(\g).$$

If $l(a)\geq 3$, we put $a=a'+\delta_i$ for a well-chosen $i$, with $l(a')\geq 2$. Then, there exists a non-zero constant $\lambda$ such that,
in $\U(\g)$, $\tdelta(v^a)=\lambda(v_i v^{a'})$. So, as $I$ is a left ideal, $v^a-\lambda v_iv^{a'} \in \g \cap I=(0)$, so $v^a=\lambda v_iv^{a'}$. 
Then $\xi_\star(v^a)=v_i*\phi(v^{a'})$. By the induction hypothesis, $\phi(v^{a'})$ belongs to $S_{\geq 2}(\g)$, left ideal for $*$, 
so $v^a$ belongs to $S_{\geq 2}(\g)$. \end{proof}\\

As a conclusion:

\begin{theo} \label{11}
Let $\g$ be a Lie algebra. There exists a prelie product on $\g$ inducing its Lie bracket, if, and only if, there exists a left ideal $I$
of $\U(\g)$ such that $\U_+(\g)=\g\oplus I$. More precisely, there exists a bijection between the sets:
\begin{itemize}
\item $\PL(\g)=\{\star$ $|$ $\star$ is a  prelie product  on $\g$ inducing the bracket of $\g\}$.
\item $\LI(\g)=\{I$ left ideal of $\U(g)$ $\mid$ $\U_+(\g)=\g \oplus I\}$.
\end{itemize}
It is given by $\Psi_\g:\LI(\g) \longrightarrow \PL(\g)$, associating to a left ideal $I$ the prelie product defined by $x\star y=\varpi_\g(xy)$, 
where $\varpi_\g$ is the canonical projection on $\g$ in the direct sum $\U_+(\g)=\g\oplus I$.
The inverse bijection is given by $\Phi_\g:\PL(\g) \longrightarrow \LI(\g)$, associating to a prelie product $\star$ the left ideal 
$\U(\g)Vect(xy-x \star y,\:x,y\in \g)$.
\end{theo}

\begin{proof} It only remains to prove that $\Phi_\g(\star)$ is generated by the elements $xy-x\star y$, $x,y \in \g$. Using the isomorphism $\xi_\star$, 
it is equivalent to prove that the left ideal $S_{\geq 2}(V)$ of $(S(V),*)$ is generated by the elements $x*y-x\star y$, $x,y\in \g$. By definition of $*$, for all $x,y\in \g$:
$$x*y-x\star y=xy+x\star y-x\star y=xy,$$
so it is equivalent to prove that the left ideal $S_{\geq 2}(V)$ is generated by $S^2(V)$. Let us denote by $J$ the left ideal generated by $S^2(V)$.
As $S_{\geq 2}(V)$ is a left ideal, $J\subseteq S_{\geq 2}(V)$. Let $v_1,\ldots,v_n \in \g$, with $n \geq 2$. Let us prove that $v_1\ldots v_n \in J$
by induction on $n$. This is obvious if $n=2$. If $n\geq 3$:
$$v_1*(v_2\ldots v_n)=v_1\ldots v_n+v_1\star (v_2\ldots v_n).$$
By lemma \ref{5}, $v_1\star (v_2\ldots v_n) \in S^{n-1}(V)$. By the induction hypothesis, as $n \geq 3$,
$v_1*(v_2\ldots v_n)$ and $v_1 \star (v_2\ldots v_n)$ belong to $J$, so $v_1\ldots v_n \in J$. As a conclusion, $S_{\geq 2}(V)=J$. \end{proof}

\begin{cor}
Let $\g$ be a Lie algebra. There exists an associative product on $\g$ inducing its Lie bracket, if, and only if, there exists a bilateral ideal $I$
of $\U(\g)$ such that $\U_+(\g)=\g\oplus I$. More precisely, there exists a bijection between the sets:
\begin{enumerate}
\item $\AS(\g)=\{\star$ $|$ $\star$ is an associative product on $\g$ inducing the bracket of $\g\}$.
\item $\BI(\g)=\{I$ bilateral ideal of $\U(g)$ $/$ $\U_+(\g)=\g \oplus I\}$.
\end{enumerate}
It is given by $\Psi_\g:\LI(\g) \longrightarrow \PL(\g)$, associating to a left ideal $I$ the associative product defined by $x\star y=\varpi_\g(xy)$, 
where $\varpi_\g$ is the canonical projection on $\g$ in the direct sum $\U_+(\g)=\g\oplus I$.
\end{cor}

\begin{proof} It is enough to verify that $\Psi_\g(\BI(\g))\subseteq \AS(\g)$ and $\Phi_\g(\AS(\g)) \subseteq \BI(\g)$.
Let $\star \in \AS(\g)$. By proposition \ref{6}, $S_{\geq 2}(\g)$ is a bilateral ideal, so is $\Phi_\g(\star)=\xi_\star^{-1}(S_{\geq 2}(\g))$.
Let $I\in \AS_{gr}(\g))$. Then $\g\approx \U_+(\g)/I$ inherits an associative product, which is $\star=\Psi_\g(I)$. \end{proof}

\subsection{Good and very good pointed modules}

\begin{defi} \label{13}
\textnormal{Let $\g$ be a Lie algebra. 
\begin{enumerate}
\item A {\it pointed $\g$-module} is a couple $(M,m)$, where $M$ is a $\g$-module and $m \in M$.
\item Let $(M,m)$ be a pointed $\g$-module. The application $\Upsilon_m$ is defined by:
$$\Upsilon_m:\left\{ \begin{array}{rcl}
\g&\longrightarrow&M\\
x&\longrightarrow&x.m.
\end{array}\right.$$
\item A pointed $\g$-module $(M,m)$ is {\it good} if the following assertions hold:
\begin{itemize}
\item $\Upsilon_m$ is injective.
\item $Im(\Upsilon_m)$ is a submodule of $M$ which does not contains $m$.
\item $M/Im(\Upsilon_m)\approx K$ (trivial $\g$-module) as a $\g$-module.
\end{itemize}
\item A pointed $\g$-module $(M,m)$ is {\it very good} if $\Upsilon_m$ is bijective.
\end{enumerate} }\end{defi}

Note that all good pointed $\g$-modules have the same dimension, so we can consider the set $\GM(\g)$ of isomorphism classes of good pointed $\g$-modules.
Similarly, we can consider the set $\VGM(\g)$ of very good pointed $\g$-modules.\\

{\bf Remark.} If $(M,m)$ is good, then it is cyclic, generated by $m$. Moreover, $M/Im(\Upsilon_m)$ is one-dimensional, trivial, generated by 
$\overline{m}=m+Im(\Upsilon_m)$.

\begin{theo} \label{14}
Let $\g$ be a Lie algebra. The following application is a bijection:
$$\Theta:\left\{\begin{array}{rcl}
\LI(\g)&\longrightarrow&\GM(\g)\\
I&\longrightarrow&(\U(\g)/I,\overline{1}).
\end{array}\right.$$
\end{theo}

\begin{proof} Let us first proove that $\Theta$ is well-defined. As $I\in \LI(\g)$, as a vector space $\U(\g)/I=(\overline{1})\oplus \g$.
For all $x,y \in \g$, $xy \in \U_+(\g)$, so $x.\overline{y} \in \U_+(\g)/I=\g$ in $\U(\g)/I$, as $I \subseteq \U_+(\g)$.
As a conclusion, the subspace $\g$ of $\U(\g)/I$ is a submodule.
Moreover, $\Upsilon_{\overline{1}}(x)=x$ for all $x \in \g$. So $\Upsilon_{\overline{1}}$ is injective, and its image is the submodule $\g$, 
so does not contain $\overline{1}$. Finally, $(\U(\g)/I)/Im(\Upsilon_{\overline{1}})\approx \U(\g)/\U_+(\g)\approx K$ as a $\g$-module.
So $(\U(\g)/I,\overline{1})$ is a good pointed $\g$-module. \\

Let us consider:
$$\Theta':\left\{\begin{array}{rcl}
\GM(\g)&\longrightarrow&\LI(\g)\\
(M,m)&\longrightarrow& Ann(m)=\{x\in \U(\g)\:\mid\:x.m=0\}.
\end{array}\right.$$
Let us prove that $\Theta'$ is well-defined. Let $(M,m)$ be a good pointed $\g$-module. Then $Ann(m)$ is a left-ideal of $\U(\g)$.
Let $x \in Ann(m)$. Then $x.m=0$, so $x.\overline{m}=0$ in $M/Im(\Upsilon_m)$, so $x.\overline{m}=\varepsilon(x)\overline{m}=0$. 
As a conclusion, $\varepsilon(x)=0$, so $Ann(m) \subseteq \U_+(\g)$. Let us now show that $\U_+(\g)=\g\oplus Ann(m)$. 
First, if $x \in \g \cap Ann(m)$, then $\Upsilon_m(x)=x.m=0$. As $\Upsilon_m$ is monic, $x=0$. If $y \in \U_+(\g)$, 
then $y.\overline{m}=\varepsilon(y)\overline{m}=0$ in $M/Im(\Upsilon_m)$, so $y.m \in Im(\Upsilon_m)$: there exists $x \in \g$, such that $y.m=x.m$. 
Hence, $y-x \in Ann(m)$, so $y=x+(y-x) \in \g+Ann(m)$. Finally, $Ann(m) \in \LI(\g)$.

Moreover, if $(M,m)$ and $(M',m')$ are isomorphic (that is to say there is an isomorphism of $\g$-modules from $M$ to $M'$ sending $m$ to $m'$),
then $Ann(m)=Ann(m')$. So $\Theta'$ is well-defined. \\

Let $(M,m)$ be a good pointed $\g$-module. Then $\Theta \circ \Theta'(M,m)=(\U(\g)/Ann(m), \overline{1}) \approx (M,m)$, as $M$ is cyclic, generated by $m$.
Let now $I\in \LI(\g)$. Then $\Theta' \circ \Theta(I)$ is the annihilator of $\overline{1}$ in $\U(\g)/I$, so is equal to $I$.
As a conclusion, $\Theta$ and $\Theta'$ are inverse bijections. \end{proof}

\begin{defi}\textnormal{The set $\GM'(\g)$ is the set of isomorphism classes of couple $((M,m),V)$, where:
\begin{itemize}
\item $(M,m)$ is a good pointed $\g$-module.
\item $V$ is a submodule of $M$ such that $M=V\oplus Im(\Upsilon_m)$.
\end{itemize}}\end{defi}

\begin{prop}
For any Lie algebra $\g$, $\GM'(\g)$ is in bijection with $\VGM(\g)$.
\end{prop}

\begin{proof}
Let $((M,m),V)$ be an element of $\GM'(\g)$. As $M/Im(\Upsilon_m)\approx K$, the complement $V$ of $Im(\Upsilon_m)$ is trivial and one-dimensional. 
As $m\notin Im(\Upsilon_m)$, $V$ admits a unique element $m'=m+x.m$, where $x \in \g$. Let us show that $(Im(\Upsilon_m), -x.m)$ is very good.
Let $m'' \in Im(\Upsilon_m)$. There exists $y \in \g$, such that $y.m=m''$. Then
$y.m'=\varepsilon(y)m'=0=y.m+y.(x.m)=m''+y.(x.m)$, so $m''=y.(-x.m)=\Upsilon_{-x.m}(y)$. Hence, $\Upsilon_{-x.m}$ is epic.

Let us assume that $\Upsilon_{-x.m}(y)=0$. Then $y.m'=\varepsilon(y)m'=0=y.m+y.(x.m)=y.m=\Upsilon_m(x)$. As $\Upsilon_m$ is monic, $y=0$,
so $\Upsilon_{-x.m}$ is also monic.

Moreover, if $((M,m),V)\approx ((M',m'),V')$, then $(Im(\Upsilon_m),-x.m)$ and $(Im(\Upsilon_{m'},-x'.m')$ are isomorphic, 
so the following application is well-defined:
$$\Lambda: \left\{\begin{array}{rcl}
\GM'(\g)&\longrightarrow& \VGM(\g)\\
((M,m),V)&\longrightarrow&(Im(\Upsilon_m),-x.m).
\end{array}\right.$$

Let now $(M,m)$ be a very good pointed $\g$-module. Let us prove that $((K\oplus M,1+m),K)$ is an element of $\GM'(\g)$.
First, for all $x \in \g$, $\Upsilon_{1+m}(x)=\varepsilon(x)1+x+m=0+\Upsilon_m(x)$, so $\Upsilon_{1+m}=\Upsilon_m$.
As a consequence, $Im(\Upsilon_{1+m})=M$ is a submodule which does not contains $1+m$, $\Upsilon_{1+m}$ is monic,
and $K$ is a complement of $Im(\Upsilon_{1+m})$. Hence, $((K\oplus M,1+m),K)\in \GM'(\g)$.
Moreover, if $(M,m) \approx (M',m')$, then $((K\oplus M,1+m),K)\approx ((K\oplus M',1+m'),K)$, so the following application is well-defined:
$$\Lambda': \left\{\begin{array}{rcl}
\VGM(\g)&\longrightarrow& \GM'(\g)\\
(M,m)&\longrightarrow&((K\oplus M,1+m),K).
\end{array}\right.$$

Let $(M,m) \in \VGM(\g)$. Then $\Lambda \circ \Lambda'((M,m))=(M,m)$. Let $((M,m),V) \in \GM'(\g)$. Then
the following application is an isomorphism of pointed $\g$-modules and sends $K$ to $V$:
$$\left\{\begin{array}{rcl}
(K\oplus Im(\Upsilon_m),1-x.m)&\longrightarrow& M\\
\lambda+m''&\longrightarrow&\lambda m'+m''.
\end{array}\right.$$
So $\Lambda$ and $\Lambda'$ are inverse bijections. \end{proof}

%% file: chap2.tex
\section{Are the complex simple Lie algebras prelie?}

The aim of this section is to prove that the complex finite-dimensional simple Lie algebras are not prelie.
In this section, the base field is the field of complex $\mathbb{C}$. We shall use the notations of \cite{Fulton} in the whole section.
In particular, if $\g$ is a simple complex finite-dimensional Lie algebra, $\Gamma_{(a_1,\ldots,a_n)}$ is the simple $\g$-module
of hightest weight $(a_1,\ldots,a_n)$.

\subsection{General results}

\begin{prop}
Let $\g$ be a finite-dimensional semi-simple Lie algebra. Then there exists a prelie product on $\g$ inducing its bracket if,
and only if, there exists a very good pointed $\g$-module.
\end{prop}

\begin{proof} As the category of finite-dimensional modules over $\g$ is semi-simple, for any good pointed $\g$-module $(M,m)$, 
there exists a $V$ such that $((M,m),V) \in \GM(\g)$. Using the bijections of the preceding section, $\PL(\g)$ is non-empty if, and only if,
$\LI(\g)$ is non-empty, if, and only if, the set $\GM(\g)$ is non-empty, if, and only if, the set $\GM'(\g)$ is non-empty, if, and only if, $\VGM(\g)$ is non-empty.
\end{proof}\\

Let $\g$ be a finite-dimensional complex simple Lie algebra. In order to prove that $\g$ is not prelie, we have to prove that $\g$ has no very good pointed modules.

\begin{lemma} \label{18}
Let $(M,m)$ be a very good $\g$-module and let $M'$ be a submodule of $M$. There exists $m' \in M$, such that the following application is epic:
$$\Upsilon_{m'}:\left\{\begin{array}{rcl}
M'&\longrightarrow&M'\\
x&\longrightarrow& x.m'.
\end{array}\right.$$
\end{lemma}

\begin{proof} As $\g$ is semi-simple, let $M''$ be a complement of $M'$ in $M$ and let $m=m'+m''$ the decomposition of $m$ in $M=M'\oplus M''$. 
Let $u \in M'$. As $\Upsilon_m$ is bijective, there exists a (unique) $x \in \g$, such that $\Upsilon_m(x)=u$. So $u+0=x.m'+xm''$ in the direct sum
$M=M'\oplus M''$. Hence, $\Upsilon_{m'}(x)=u$ and $\Upsilon_{m'}$ is epic. \end{proof}

\begin{cor} \label{19}
Let $(M,m)$ be a very good $\g$-module. Then $M$ is not isomorphic to $(\g,ad)$ and its trivial component is $(0)$.
\end{cor}

\begin{proof} First, note that for all $x \in \g$, $((\g,ad),x)$ is not very good: it is obvious if $x=0$ and if $x \neq 0$, $\Upsilon_x(x)=[x,x]=0$,
so $\Upsilon_x$ is not monic. 

If $(M,m)$ is very good, let us consider its trivial component $M'$. By the preceding lemma, there exists $m' \in M'$,
such that $\Upsilon_{m'}:M'\longrightarrow M'$ is epic. As $M'$ is trivial, $\Upsilon_{m'}=0$, and finally $M'=(0)$.  \end{proof}

\subsection{Representations of small dimension}

We first give the representations of $\g$ with dimension smaller than the dimension of $\g$. 
We shall say that a $\g$-module $M$ is {\it small} if:
\begin{itemize}
\item $M$ is simple.
\item The dimension of $M$ is smaller than the dimension of $\g$.
\item $M$ is not trivial and not isomorphic to $(\g,ad)$.
\end{itemize}

{\bf Remark.} By corollary \ref{19}, any very good module is the direct sum of small modules. \\

Let $n$ be the rank of $\g$ and let $\Gamma_{(a_1,\ldots,a_n)}$ be the simple $\g$-module of highest weight $(a_1,\ldots,a_n)$.
The dimension of $\Gamma_{(a_1,\ldots,a_n)}$ is given in \cite{Fulton} in chapter 15 for $\mathfrak{sl}_{n+1}$, (formula 15.17),
in chapter 24 for $\mathfrak{sp}_{2n}$, $\mathfrak{so}_{2n}$ and $\mathfrak{so}_{2n+1}$ (exercises 24.20, 24.30 and 24.42).
It turns out from these formulas that if $b_1\leq a_1, \ldots b_n \leq a_n$, then $dim(\Gamma_{(b_1,\ldots,b_n)})\leq dim(\Gamma_{(a_1,\ldots,a_n)})$,
with equality if, and only if, $(b_1,\ldots,b_n)=(a_1,\ldots,a_n)$. Direct computations of $dim(\Gamma_{(0,\ldots,0,3,0,\ldots,0)})$ 
then shows that it is enough to compute the dimensions of $\Gamma_{(a_1,\ldots,a_n)}$ with all $a_1\leq 2$. With the help of a computer, we obtain:

$$\begin{array}{|c|c|c|c|}
\hline \g&dim(\g)&\mbox{small modules $V$}&dim(V)\\
\hline \hline \mathfrak{sl}_n,\:n\geq 9&n^2-1&\Gamma_{(1,0,\ldots,0)},\:\Gamma_{(0,\ldots,0,1)}&n\\
&&\Gamma_{(0,1,0,\ldots,0)},\:\Gamma_{(0,\ldots,0,1,0)}&n(n-1)/2\\
&&\Gamma_{(2,0,\ldots,0)},\:\Gamma_{(0,\ldots,0,2)}&(n+1)/2\\
\hline \mathfrak{sl}_n,\:3\leq n \leq 8&n^2-1&\Gamma_{(0,\ldots,0,1,\ldots,0)}, \mbox{the $1$ in position $i$}&\binom{n}{i}\\
&&\Gamma_{(2,0,\ldots,0)},\:\Gamma_{(0,\ldots,0,2)}&(n+1)/2\\
\hline \mathfrak{sl}_2&3&\Gamma_1&2\\
\hline \hline \mathfrak{sp}_{2n}, n=2 \mbox{ or }n \geq 4&n{2n+1}&\Gamma_{(1,0,\ldots,0)}&2n\\
&&\Gamma_{(0,1,0,\ldots,0)}&{2n+1}(n-1)\\
\hline \mathfrak{sp}_6&21&\Gamma_{(1,0,0)}&6\\
&&\Gamma_{(0,1,0)},\:\Gamma_{(0,0,1)}&14\\
\hline \hline \mathfrak{so}_{2n},\:n\geq 8&n(2n-1)&\Gamma_{(1,0,\ldots,0)}&2n\\
\hline \mathfrak{so}_{2n},\:4\leq n\leq 7&n(2n-1)&\Gamma_{(1,0,\ldots,0)}&2n\\
&&\Gamma_{(0,0,1,0,\ldots,0)},\:\Gamma_{(0,\ldots,0,1)}&2^{n-1}\\
\hline \mathfrak{so}_6&15&\Gamma_{(1,0,0)}&6\\
&&\Gamma_{(0,1,0)},\:\Gamma_{(0,0,1)}&4\\
&&\Gamma_{(0,2,0)},\:\Gamma_{(0,0,2)}&10\\
\hline \hline \mathfrak{so}_{2n+1},\:n\geq 7&n{2n+1}&\Gamma_{(1,0,\ldots,0)}&2n+1\\
\hline \mathfrak{so}_{2n+1},\:2\leq n\leq 6&n{2n+1}&\Gamma_{(1,0,\ldots,0)}&2n+1\\
&&\Gamma_{(0,\ldots,0,1)}&2^n\\
\hline \hline \mathfrak{e}_6&78&\Gamma_{(1,0,0,0,0,0)}&27\\
\hline \mathfrak{e}_7&133&\Gamma_{(1,0,0,0,0,0,0)}&56\\
\hline \mathfrak{e}_8&78&\times&\times\\
\hline \mathfrak{f}_4&52&\Gamma_{(1,0,0,0)}&26\\
\hline \mathfrak{g}_2&14&\Gamma_{(1,0)}&7\\
\hline  \end{array}$$

\subsection{Proof in the generic cases}

\begin{prop}
Let $\g$ be $\mathfrak{sl}_n$, with $n \neq 6$, or $\mathfrak{so}_{2n}$, with $n \geq 2$, or $\mathfrak{sp}_{2n}$, with $n \geq 2$,
or $\mathfrak{e}_n$, with $6\leq n\leq 8$. Then $\g$ is not prelie.
\end{prop}

\begin{proof} Let $\g$ be a simple finite-dimensional Lie algebra, with a prelie product inducing its Lie bracket.
From the preceding results, $\g$ has a very good pointed $\g$-module $(M,m)$, and the submodules appearing in the decomposition of $M$ 
into simples modules are all small. Let us show that it is not possible in these different cases.\\

{\it First case.} If $\g=\mathfrak{sl}_n$, with $n\geq 9$, then the decomposition of a very good pointed $\g$-module into simples would have $a$ submodules 
of dimension $n$, $b$ submodules of dimensions $n(n-1)/2$ and $c$ submodules of dimension $n(n+1)/2$. So:
$$n^2-1=an+b\frac{n(n-1)}{2}+c \frac{n(n+1)}{2}.$$
Let us assume that $n$ has a prime factor $p \neq 2$. Then $p \mid n$, $p\mid n(n\pm 1)$ and $2 \mid n(n \pm 1)$,
so $p \mid \frac{n(n\pm 1)}{2}$. As a consequence, $p\mid n^2-1$, so $p\mid 1$: this is absurd. So $n=2^k$ for a certain $k \geq 4$, as $n \geq 9$.
Then:
$$a2^k+b2^{k-1}(2^k-1)+c2^{k-1}(2^k+1)=2^{2k}-1,$$
so $2^{k-1}\mid 2^{2k}-1$ and $k=1$: contradiction. \\

{\it Second case.} Similarly, if $\g=\mathfrak{sp}_{2n}$ with $n \neq 3$, we would have $a,b\in \mathbb{N}$ such that $a2n+b{2n+1}(n-1)=n{2n+1}$. 
If $b\geq 2$, then $n{2n+1}\geq 2{2n+1}(n-1)$, so ${2n+1}(n-2) \leq 0$ and $n \leq 2$: contradiction. So $b=0$ or $1$.
If $b=0$, then $a=n+\frac{1}{2}$: absurd. If $b=1$, then $a=1+\frac{1}{2n}$: absurd. \\

{\it Third case.} If $\g=\mathfrak{so}_{2n}$, with $n \geq 8$, then we would have $2n \mid n(2n-1)$, so $n-\frac{1}{2}\in \mathbb{N}$: absurd.\\

{\it Other cases.} Let us sum up the possible dimensions of the small modules in an array:
$$\begin{array}{|c|c|c|}
\hline \g&\mbox{dimensions of the small modules}&dim(\g)\\
\hline \mathfrak{sl}_2&2&3\\
\hline \mathfrak{sl}_3&3,6&8\\
\hline \mathfrak{sl}_4&4,6,10&15\\
\hline \mathfrak{sl}_5&5,10,15&24\\
\hline \mathfrak{sl}_7&7,21,28,35&48\\
\hline \mathfrak{sl}_8&8,28,36,56&63\\
\hline \mathfrak{so}_6&4,6,10&15\\
\hline \mathfrak{so}_8&8&28\\
\hline \mathfrak{so}_{10}&10,16&45\\
\hline \mathfrak{so}_{12}&12,32&66\\
\hline \mathfrak{so}_{14}&14,64&91\\
\hline \mathfrak{sp}_6&6,14&21\\
\hline \mathfrak{e}_6&27&78\\
\hline \mathfrak{e}_7&56&133\\
\hline \mathfrak{e}_8&\times&248\\
\hline \end{array}$$
In all cases except $\mathfrak{so}_6$, (there is a prime integer which divides all the dimensions of the small submodules but not the dimension of $\g$),
or (there is a unique dimension of small modules and it does not divide the dimension of $\g$), or (there is no small modules). So $\g$ is not prelie.
For $\mathfrak{so}_{12}$, it comes from the fact that $4$ divides $12$ and $32$ but divides not $66$.\\

So none of these Lie algebras is prelie.  \end{proof}

\subsection{Case of $\mathfrak{so}_{2n+1}$}

\label{s2.5}

We assume that $n \geq 2$. Let us recall that:
$$\mathfrak{so}_{2n+1}=\left\{
\left(\begin{array}{ccc}
A&B&-^tF\\
C&-^tA&-^tE\\
E&F&0
\end{array} \right) \:\mid\:B,C \mbox{ skew symmetric}\right\},$$
 where $A,B,C$ are $n\times n$ matrices, $E,F$ are $1 \times n$ matrices.

\begin{lemma}
Let $(M,m)$ be a very good pointed module over $\g=\mathfrak{so}_{2n+1}$. Then $M$ is isomorphic to
$M_{2n+1,n}(\mathbb{C})$ as a $\g$-module, with the action of $\g$ given by the (left) matricial product.
\end{lemma}

\begin{proof}
Let us first assume that $n\geq 7$. Then $\g$ has a unique small module of dimension $2n+1$, which is the standard representation $\mathbb{C}^{2n+1}$.
So if $M$ is a direct sum of copies of this module; comparing the dimension, there are necessarily $n$ copies of this module, so
$M\approx M_{2n+1,n}(\mathbb{C})$.

If $n \leq 6$, let us sum up the possible dimensions of the small modules in an array:
$$\begin{array}{|c|c|c|}
\hline \g&\mbox{dimensions of the small modules}&dim(\g)\\
\hline \mathfrak{so}_5&5,4&10\\
\hline \mathfrak{so}_7&7,8&21\\
\hline \mathfrak{so}_9&9,16&36\\
\hline \mathfrak{so}_{11}&11,32&55\\
\hline \mathfrak{so}_{13}&13,64&78\\
\hline \end{array}$$
In all these cases, the only possible decomposition of $M$ is $n$ copies of the standard representation of dimension $2n+1$, so the conclusion also holds.
\end{proof}\\

The elements of $M_{2n+1,n}(\mathbb{C})$ will be written as $\left(\substack{X\\Y\\z}\right)$, where $X,Y \in M_n(\mathbb{C})$ and $z \in M_{1,n}(\mathbb{C})$.

\begin{prop}
For any $m=\left(\substack{X\\Y\\z}\right)\in M_{2n+1,n}(\mathbb{C})$, $(M_{2n+1,n}(\mathbb{C}),m)$ is not very good.
\end{prop}

\begin{proof} Let us denote:
$$\M=\{A \in M_{2n+1,n}((\mathbb{C})\:\mid \:(M_{2n+1,n}(\mathbb{C}),A)\mbox{ is very good}\}.$$
Let us recall that:
$$SO_{2n+1}=\left\{B \in M_{2n+1,2n+1}(\mathbb{C})\:\mid \: ^tB \left(\begin{array}{ccc}
0&I&0\\I&0&0\\0&0&I \end{array}\right)B=\left(\begin{array}{ccc}
0&I&0\\I&0&0\\0&0&I \end{array}\right)\right\}.$$
It is well-known that if $B \in SO_{2n+1}$ and $A \in \mathfrak{so}_{2n+1}$, then $BAB^{-1} \in \mathfrak{so}_{2n+1}$. \\

{\it First step.} Let us prove that $m\in \M$ if, and only if, $Bm \in \M$ for all $B \in SO_{2n+1}$.
Let us assume that $m\in \M$ and let $B \in SO_{2n+1}$. Let us take $A \in Ker(\Upsilon_{Bm})$. Then $ABm=0$,
so $B^{-1}ABm=0$ and $B^{-1}AB \in Ker(\Upsilon_{m})$. As $m\in\M$, $\Upsilon_m$ is bijective, so $A=0$: $\Upsilon_{Bm}$ is monic.
As $M$ and $\mathfrak{so}_{2n+1}$ have the same dimension, $\Upsilon_{Bm}$ is bijective.\\

{\it Second step.} Let us prove that $m\in \M$ if, and only if, $mP\in \M$ for all $P \in GL(n)$. Indeed, $m\longrightarrow mP$
is an isomorphism of $\g$-modules for all $P \in GL(n)$. \\

{\it Third step}. Let us prove that $\left(\substack{X\\Y\\z}\right)\in \M$ if, and only if, 
$\left(\substack{QXP\\^tQ^{-1}YP\\zP}\right)\in \M$ for all $P,Q \in GL(n)$.
Indeed, if $\left(\substack{X\\Y\\z}\right)\in \M$ and if $P,Q \in GL(n)$, by the first and second steps:
$$\left(\begin{array}{rcl}
Q&0&0\\
0&^tQ^{-1}&0\\
0&0&1
\end{array}\right)\left(\begin{array}{c}
X\\Y\\z\end{array}\right)P\in \M,$$
as the left-multiplying matrix is an element of $SO_{2n+1}$. \\

{\it Fourth step}. Let $\left(\substack{X\\Y\\z}\right)\in \M$, let us prove that $X$ is invertible. If $X$ is not invertible,
for a good choice of $P$ and $Q$, $QXP$ has its first row and first column equal to $0$. By the third step,
$m'=\left(\substack{QXP\\^tQ^{-1}YP\\zP}\right)\in \M$, but $\Upsilon_{m'}(E_{n+1,2}-E_{n+2,1})=0$, where $E_{i,j}$ is the elementary matrix
with only a $1$ in position $(i,j)$: this is a contradiction. So $X$ is invertible.\\

{\it Last step.} Let us assume that $\M$ is not empty and let $m=\left(\substack{X\\Y\\z}\right)\in \M$. Then $X$ is invertible.
For a good choice of $P$ and $Q$, we obtain an element $m'=\left(\substack{I_n\\Y'\\z'}\right)\in \M$. Let us then choose a non-zero skew-symmetric 
matrix $B \in M_n(\mathbb{C})$ (this exists as $n \geq 2$), then the following element is in $\mathfrak{so}_{2n+1}$:
$$A=\left(\begin{array}{rcl}
-BY&B&0\\
^tYBY&-^tYB&0\\
0&0&0
\end{array}\right)\in \mathfrak{so}_{2n+1}.$$
An easy computation shows that $A.m'=0$, so $A \in Ker(\Upsilon_{m'})$: contradiction, $m'\notin \M$. So $\M$ is empty. \end{proof}

\begin{cor}
For all $n \geq 2$, $\mathfrak{so}_{2n+1}$ is not prelie.
\end{cor}

\subsection{Proof for $\mathfrak{sl}_6$ and $\mathfrak{g}_2$}

The Lie algebra $\mathfrak{sl}_6$ has dimension $35$, and the possible dimensions of its small modules are $6$, $15$, $20$.
So if $(M,m)$ is a very good pointed module, $M$ is the direct sum of a simple of dimension $15$ and a simple of dimension $20$. 
As there is only one simple of dimension $20$, that is to say $\Lambda^3(V)$ where $V$ is the standard representation, $M$ contains 
$\Lambda^3(V)$. By lemma \ref{18}, in order to prove that $\mathfrak{sl}_6$ is not prelie, it is enough to prove that for any $m\in \Lambda^3(V)$,
$\Upsilon_m$ is not epic. This essentially consists to show that the rank of a certain $20\times 35$ matrix is not $20$ 
and this can be done directly using MuPAD pro 4, see section \ref{s4.2}. \\

The proof for $\mathfrak{g}_2$ is similar: if $(M,m)$ is a very good module, then $M\approx V \oplus V$, where $V$ is the only small module of $\g$,
that is to say its standard representation. So $M\approx M_{7,2}(\mathbb{C})$ as a $\g$-module. It remains to show that for any $m$,
$(M_{7,2}(\mathbb{C}),m)$ is not very good. This essentially consists to show that a certain $14\times 14$ matrix is not invertible 
and this can be done directly using MuPAD pro 4, see section \ref{s4.3}.\\

The proof for $\mathfrak{f}_4$ would be similar: if $(M,m)$ is a very good module, then $M\approx V \oplus V$, where $V$ is the only small module of $\g$,
that is to say its standard representation. So $M\approx M_{26,2}(\mathbb{C})$ as a $\mathfrak{f}_4$-module. It would remain to show that for any $m$,
$(M_{26,2}(\mathbb{C}),m)$ is not very good. This would consist to show that a certain $52 \times 52$ matrix is not invertible.\\

We finally prove:
\begin{theo} 
Let $\g$ be a simple, finite-dimensional complex Lie algebra. If it is not isomorphic to $\mathfrak{f}_4$, it is not prelie.
\end{theo}

We conjecture:
\begin{theo} 
Let $\g$ be a simple, finite-dimensional complex Lie algebra. Then it is not prelie.
\end{theo}

{\bf Remark.} As a corollary, we obtain that $\g$ is not associative. This result is also proved in a different way in \cite{Kubo}, with the help of
compatible products. Indeed, we have the following equivalences:

\begin{tabular}{rcl}
&&the prelie product $\star$ is compatible\\
&$\Longleftrightarrow$& $\forall x,y,z \in \g$, $[x,y\star z]=[x,y]\star z+y\star[x,z]$\\
&$\Longleftrightarrow$& $\forall x,y,z \in \g$, $x\star(y\star z)-(y\star z)\star x-(x\star y)\star z+(y\star x) \star z-y\star(x\star z)+y \star(z \star x)$\\
&$\Longleftrightarrow$& $\forall x,y,z \in \g$, $-(y\star z)\star x-+y \star(z \star x)=0$\\
&$\Longleftrightarrow$& $\star$ is associative.
\end{tabular}

As from \cite{Kubo}, the only admissible associative product on $\g$ is $0$, $\g$ is not associative.

%% file: chap3.tex
\section{Dendriform products on cofree coalgebra}

\subsection{Preliminaries and results on tensor coalgebras}

Let $V$ be a vector space. The tensor algebra $T(V)$ has a coassociative coproduct, given for all $v_1,\ldots,v_n \in V$ by:
$$\Delta(v_1\ldots v_n)=\sum_{i=0}^n v_1\ldots v_i \otimes v_{i+1}\ldots v_n.$$
Let us recall the following facts:
\begin{enumerate}
\item Let us fix a basis $(v_i)_{i\in I}$ of $V$. We define $\TI$ as the set of words $w$ in letters the elements of $I$.
For an element $w=i_1\ldots i_k \in \SI$, we put $l(w)=k$, and:
$$v_w=\prod_{i\in I} v_{i_1}\ldots v_{i_k}.$$
Then $(v_w)_{w\in \TI}$ is a basis of $T(V)$, and the coproduct is given by:
$$\Delta(v_w)=\sum_{w_1w_2=w}v_{w_1} \otimes v_{w_2}.$$
\item Let $T_+(V)$ be the augmentation ideal of $T(V)$. It is given a coassociative, non counitary coproduct $\tdelta$ defined by
$\tdelta(x)=\Delta(x)-x\otimes 1-1\otimes x$ for all $x \in T_+(V)$. In other terms, putting $\TIp=\TI-\{\emptyset\}$, 
$(v_w)_{w\in \TIp}$ is a basis of $T_+(V)$ and:
$$\tdelta(v_w)=\sum_{\substack{w_1w_2=w\\w_1,w_2 \in \TIp}} v_{w_1} \otimes v_{w_2}.$$
\item Let $\tdelta^{(n)}:T_+(V) \longrightarrow T_+(V)^{\otimes (n+1)}$ be the $n$-th iterated coproduct of $T_+(V)$. Then:
$$Ker\left(\tdelta^{(n)}\right)=\bigoplus_{k=1}^n T^k(V).$$
In particular, $Prim(T(V))=V$.
\item $(T^n(V))_{n\in \mathbb{N}}$ is a gradation of the coalgebra $T(V)$.
\end{enumerate}

\begin{lemma} 
In $T_+(V)\otimes T_+(V)$, $Ker(\tdelta \otimes Id-Id \otimes \tdelta)=Im(\tdelta)$.
\end{lemma}

\begin{proof} $\supseteq$ comes from the coassociativity of $\tdelta$.

$\subseteq$: let us consider $\displaystyle X=\sum_{w_1,w_2\in \TIp} x_{w_1,w_2} v_{w_1}\otimes v_{w_2} \in 
Ker(\tdelta \otimes Id-Id \otimes \tdelta)=Im(\tdelta)$. Then:
\begin{eqnarray*}
(\tdelta \otimes Id) \circ \tdelta(X)&=&\sum_{w_1,w_2,w_3 \in \TIp} x_{w_1w_2,w_3} v_{w_1} \otimes v_{w_2} \otimes v_{w_3}\\
=(Id \otimes \tdelta) \circ \tdelta(X)&=&\sum_{w_1,w_2,w_3 \in \TIp} x_{w_1,w_2w_3} v_{w_1} \otimes v_{w_2} \otimes v_{w_3}.
\end{eqnarray*}
So, for all $w_1,w_3 \in \TIp$, $w_2\in \TI$, $x_{w_1w_2,w_3}=x_{w_1,w_2w_3}$.
In particular, $x_{i_1\ldots i_k,j_1\ldots j_k}=x_{i_1,i_2\ldots i_kj_1\ldots j_l}$ for all $i_1,\ldots,i_k,j_1,\ldots,j_l \in I$.
We denote by $x_{i_1\ldots i_kj_1\ldots j_l}$ this common value. Then:
$$\Delta(X)=\sum_{w\in \TIp} x_{w} \sum_{\substack{w_1w_2=w\\w_1,w_2 \in \SI_+}} v_{w_1} \otimes v_{w_2}=\sum_{w\in \TIp} x_{w} \tdelta(x_w).$$
So $X \in Im(\tdelta)$. \end{proof}

\begin{lemma}
Let $W$ be a subspace of $T_+(V)$, such that $T(V)=(1)\oplus V \oplus W$. There exists a coalgebra endomorphism $\phi$ of $T(V)$,
such that $\phi_{\mid V}=Id_V$ and $\displaystyle \phi\left(\bigoplus_{n\geq 2} T^n(V)\right)=W$. Moreover, $\phi$ is an automorphism.
\end{lemma}

\begin{proof} Similar to the proof of lemma \ref{3}. \end{proof}

\begin{cor} \label{28}
Let $C$ be a cofree coalgebra and let $W\subset C_+$, such that $C=(1)\oplus Prim(C)\oplus W$. There exists a unique coalgebra isomorphism
$\phi:T(Prim(C))\longrightarrow C$ such that $\phi_{\mid Prim(C)}=Id_{Prim(C)}$ and $\phi(T_{\geq 2}(Prim(C))=W$.
\end{cor}

\subsection{Left ideal associated to a dendriform Hopf algebra}

Let $A$ be a dendriform Hopf algebra \cite{Loday1,Loday3}, that is to say the product (denoted by $*$) of $A$ can be split on $A_+$ as $*=\prec+\succ$, with:
\begin{enumerate}
\item For all $x,y,z \in A_+$:
\begin{eqnarray}
\label{E1} (x \prec y) \prec z&=&x \prec (y*z),\\
\label{E2} (x \succ y) \prec z&=&x \succ(y\prec z),\\
\label{E3} (x*y) \succ z&=&(x \succ y) \succ z.
\end{eqnarray}
\item For all $a,b \in A_+$:
\begin{eqnarray}
\label{E4} \tdelta(a \prec b)&=&a'*b'\otimes a''\prec b''+a'*b\otimes a''+b' \otimes a \prec b''+a' \otimes a''\prec b+b \otimes a,\\
\label{E5} \tdelta(a \succ b)&=&a'*b' \otimes a'' \succ b''+a*b' \otimes b''+b'\otimes a \succ b''+a' \otimes a''\succ b+a \otimes b.
\end{eqnarray}
We used the following notations: $A_+$ is the augmentation ideal of $A$ and $\tdelta:A_+\longrightarrow A_+\otimes A_+$ is the 
coassociative coproduct defined by $\tdelta(a)=\Delta(a)-a\otimes 1-1\otimes a$ for all $a \in A_+$.
\end{enumerate}

\begin{prop} \label{29}
Let $A$ be a dendriform Hopf algebra. Then the following application is a monomorphism of coalgebras:
$$\Theta_A:\left\{\begin{array}{rcl}
T(Prim(A))&\longrightarrow&A\\
v_1\ldots v_n &\longrightarrow & v_n \prec(v_{n-1} \prec(\ldots \prec(v_2 \prec v_1)\ldots )
\end{array}\right.$$
Moreover, if $A$ is connected as a coalgebra, then $\Theta_A$ is an isomorphism, so $A$ is a cofree coalgebra.
\end{prop}

\begin{proof} For all $v_1,\ldots,v_n \in Prim(A)$, we put:
$$\omega(v_1,\ldots,v_n)=v_n \prec(v_{n-1} \prec(\ldots \prec(v_2 \prec v_1)\ldots ).$$
An easy induction using (\ref{E4}) proves that:
$$\Delta(\omega(v_1,\ldots,v_n))=\sum_{i=0}^n \omega(v_1,\ldots,v_i) \otimes \omega(v_{i+1},\ldots,v_n).$$
So $\Theta_A$ is a morphism of coalgebras. Let us prove that this morphism is injective. If not, its kernel would contain primitive elements of $T(Prim(A))$, 
that is to say elements of $Prim(A)$: absurd. If $A$ is connected, this morphism is surjective: let us take $x \in A$, let us prove that $x \in Im(\Theta)$. 
As $A$ is connected, for all $x \in A$,  there exists $n \geq 1$ such that $\tdelta^{(n)}(x)=0$. Let us proceed by induction on $n$. If $n=1$, 
then $x \in \g \subseteq Im(\Theta)$.  If $n \geq 2$, then $\tdelta^{(n-1)}(x)\in \g^{\otimes n}$. Let us put:
$$\tdelta^{(n-1)}(x)=\sum_{i_1,\ldots,i_n \in I} a_{i_1\ldots i_n} v_{i_1}\otimes \ldots \otimes v_{i_n}.$$
Then:
$$\tdelta^{(n-1)}(x)=\tdelta^{(n-1)}\left(\sum_{i_1,\ldots,i_n \in I} a_{i_1\ldots i_n} v_{i_1\ldots i_n}\right).$$
By the induction hypothesis, $\displaystyle x-\sum_{i_1,\ldots,i_n \in I} a_{i_1\ldots i_n} v_{i_1\ldots i_n}\in Im(\Theta)$, so $x \in Im(\Theta)$.
As a conclusion, $\Theta$  is an isomorphism, so $(v_w)_{w \in \TI}$ is a basis of $A$. \end{proof}

\begin{prop} \label{30}
Let $A$ be a dendriform Hopf algebra, connected as a coalgebra. Let us put $A^{\prec 2}=Prim(A) \prec A_+$. Then:
\begin{enumerate}
\item $A_+=Prim(A) \oplus A^{\prec 2}$.
\item $A^{\prec 2}=A_+ \prec A_+$.
\item $A^{\prec 2}$ is a left ideal of $A$.
\end{enumerate}
\end{prop}

\begin{proof} \begin{enumerate}
\item Using $\Theta_A$:
$$\Theta_A(Prim(A))=Prim(A),\:\Theta_A\left(\bigoplus_{n\geq 2} Prim(A)^{\otimes n}\right)=A^{\prec 2},\:
\Theta_A\left(\bigoplus_{n\geq 1} Prim(A)^{\otimes n}\right)=A_+.$$
As $\Theta_A$ is an isomorphism, we obtain the result.

\item $\subseteq$. As $Prim(A)\subseteq A_+$, $A^{\prec 2} \subseteq A_+\prec A_+$. 

$\supseteq$. Note that $A_+=Vect\left(\omega(v_1,\ldots,v_n)\:\mid\: n \geq 1,v_1,\ldots,v_n \in Prim(A)\right)$.
Let $x=\omega(v_1,\ldots,v_m)$ and $y=\omega(w_1,\ldots,w_n) \in A_+$. We put $x'=\omega(v_1,\ldots,v_{m-1})$. Then, by (\ref{E1}):
$$x \prec y=(v_m \prec x') \prec y=v_m \prec (x'*y).$$
As $x'*y \in A_+$, $x \prec y \in A^{\prec 2}$.

\item It is enough to prove that $A_+ *A^{\prec 2} \subseteq A^{\prec 2}$. Let $x,y,z \in A_+$.
$$x*(y\prec z)=x \prec(y\prec z)+(x \succ y)\prec z\in A_+\prec A_+=A^{\prec 2}.$$
So $A^{\prec 2}$ is a left ideal of $A$. 
\end{enumerate} \end{proof}

\subsection{Dendriform products on a tensorial coalgebra}

Let us now consider an associative product $*$ on the coalgebra $T(V)$, such that:
\begin{enumerate}
\item $(T(V),*,\Delta)$ is a Hopf algebra.
\item $\displaystyle T(V)_{\geq 2}=\bigoplus_{n \geq 2}T^n(V)$ is a left ideal of $(T(V),*)$.
\end{enumerate}

\begin{prop} \label{31}
We define a product $\prec$ on $T_+(V)$ in the following way: for all $v,v_1,\ldots,v_n \in V$,  $w \in T_+(V)$,
$$ \left\{ \begin{array}{rcl}
v \prec w&=&wv,\\
(v_1\ldots v_n) \prec w&=&((v_1\ldots v_{n-1})*w)v_n.
\end{array}\right.$$
We also put $\succ=*-\prec$. Then $(T(V),\prec,\succ,\tdelta)$ is a dendriform Hopf algebra.
\end{prop}

\begin{proof}  Let us first prove (\ref{E1}). Let $u_1,\ldots,u_k,v_1,\ldots,v_l, w_1,\ldots,w_m \in V$.
\begin{eqnarray*}
(u_1\ldots u_k \prec v_1\ldots v_l) \prec w_1\ldots w_m&=&((u_1\ldots u_{k-1}*v_1\ldots v_l)u_k) \prec w_1\ldots w_m\\
&=&(u_1\ldots u_{k-1}*v_1\ldots v_l *w_1\ldots w_m)u_k\\
&=&u_1\ldots u_k \prec ( v_1\ldots v_l* w_1\ldots w_m).
\end{eqnarray*}
Let us now prove (\ref{E4}). We take $a=u_1\ldots u_k$, $b=v_1\ldots v_l \in T_+(V)$. We denote $\tilde{a}=u_1\ldots u_{k-1}$.
\begin{eqnarray*}
\tdelta(a \prec b)&=&\tdelta((\tilde{a}*b)u_k)\\
&=&\tilde{a}*b \otimes u_k+\tdelta(\tilde{a}*b)(1 \otimes u_k)\\
&=&\tilde{a}*b \otimes u_k+\tilde{a} \otimes b u_k+b \otimes \tilde{a} u_k+\tilde{a}'*b \otimes \tilde{a}''u_k\\
&&+\tilde{a}' \otimes (\tilde{a}''*b)u_k+\tilde{a}*b'\otimes b''u_k+b'\otimes (\tilde{a}*b'')u_k+\tilde{a}'*b'\otimes (\tilde{a}''*b'')u_k.
\end{eqnarray*}
Moreover, using (\ref{E1}):
\begin{eqnarray*}
a'*b' \otimes a''\prec b''&=&\tilde{a}*b' \otimes b''u_k+\tilde{a}'*b' \otimes (\tilde{a}''*b'')u_k,\\
a'*b\otimes a''&=&\tilde{a}*b\otimes u_k+\tilde{a}'*b\otimes \tilde{a}''u_k,\\
b' \otimes a \prec b''&=&b' \otimes (\tilde{a}*b'')u_k,\\
a' \otimes a''\prec b&=&\tilde{a} \otimes bu_k+\tilde{a}' \otimes (\tilde{a}''*b)u_k,\\
b\otimes a&=&b \otimes \tilde{a}u_k.
\end{eqnarray*}
So (\ref{E4}) is satisfied. As $T(V)$ is a Hopf algebra, (\ref{E4})+(\ref{E5}) is satisfied, so (\ref{E5}) also is.

Let us prove (\ref{E2}). For all $x,y,z \in T_+(V)$, we put $\phi(x,y,z)=(x \succ y)\prec z-x \succ(y\prec z)$. A direct computation 
using (\ref{E4}) and (\ref{E5}) shows that:
\begin{eqnarray*}
\tdelta(\phi(x,y,z))&=&x'y'z'\otimes \Phi(x'',y'',z'')+y'z'\otimes \Phi(x,y'',z'')+x'z'\otimes \Phi(x'',y,z'')\\
&&+x'y'\otimes \Phi(x'',y'',z)+z'\otimes \Phi(x,y,z'')+y'\otimes \Phi(x,y'',z)+x'\otimes \Phi(x'',y,z).
\end{eqnarray*}
Moreover:
$$\phi(x,y,z)=(x\succ y)\prec z-x*(y\prec z)+x \prec (y\prec z).$$
By definition of $\prec$, $(x \succ y) \prec z$ and $x \prec (y\prec z)\in T_{\geq 2}(V)$. In the same way, $y\prec z \in T_{\geq 2}(V)$, 
left ideal of $T(V)$, so $x*(y\prec z) \in T_{\geq 2}(V)$. finally, $\phi(x,y,z)\in T_{\geq 2}(V)$.

Let us now prove that $\phi(x,y,z)=0$ by induction on $n=l(x)+l(y)+l(z)$. If $n=3$, then $x,y,z\in V$, so are primitive.
So $\tdelta(\phi(x,y,z))=0$, and $\phi(x,y,z)\in V$. By the preceding point, $\phi(x,y,z)\in V\cap T_{\geq 2}(V)=(0)$, so $\phi(x,y,z)=0$. 
Let us assume the result at all rank $<n$. By the induction hypothesis applied to $x',y',z'$ and others, $\tdelta(\phi(x,y,z))=0$, 
so $\phi(x,y,z)\in V\cap T_{\geq 2}(V)=(0)$. Hence, (\ref{E2}) is satisfied. As $*$ is associative, (\ref{E1})+(\ref{E2})+(\ref{E3}) is satisfied, 
so (\ref{E3}) also is. \end{proof} \\

{\bf Remark.} In the dendriform Hopf algebra $T(V)$, $T(V)^{\prec 2}=T(V)_{\geq 2}$.\\

By the dendriform Cartier-Quillen-Milnor-Moore theorem, $T(V)$ is now the dendriform enveloping algebra of the brace algebra $V=Prim(T(V))$. 
By \cite{Ronco}, the brace structure on $V$ induced by the dendriform structure of $A$ is given, for all $a_1,\ldots,a_n \in V$, by:
\begin{eqnarray*}
&&\langle a_1,\ldots,a_n \rangle\\
&=&\sum_{i=0}^{n-1}(-1)^{n-1-i} (a_1\prec(a_2\prec(\ldots \prec a_i)\ldots)\succ a_n \prec(\ldots (a_{i+1}\succ a_{i+2})\succ \ldots )\succ a_{n-1})\\
&=&\sum_{i=0}^{n-1}(-1)^{n-1-i} (a_1\ldots a_i)\succ a_n \prec(\ldots (a_{i+1}\succ a_{i+2})\succ \ldots )\succ a_{n-1})\\
&=&\pi_V\left(\sum_{i=0}^{n-1}(-1)^{n-1-i} (a_1\ldots a_i)\succ a_n \prec(\ldots (a_{i+1}\succ a_{i+2})\succ \ldots )\succ a_{n-1})\right),
\end{eqnarray*}
where we denote by $\pi_V$ the canonical projection on $V$ in $T(V)$. We obtain,as $T(V)^{\prec 2}=Ker(\pi_V)$:
\begin{eqnarray*}
\langle a_1,\ldots,a_n \rangle&=&\pi_V((a_1\ldots a_{n-1})\succ a_n)\\
&=&\pi_V((a_1\ldots a_{n-1})*a_n)-\pi_V((a_1\ldots a_{n-1})\prec a_n)\\
&=&\pi_V((a_1\ldots a_{n-1})*a_n).
\end{eqnarray*}
In other terms, identifying $V$ and $T_+(V)/T(V)_{\geq 2}$, $V$ becomes a left $(T(V),*)$-module, and the brace structure of $V$ is given
by this module structure.

\subsection{Dendriform structures on a cofree coalgebra}

\begin{theo} \label{32}
Let $A$ be a cofree coalgebra. We define:
\begin{enumerate}
\item $\DD(A)=\{(\prec,\succ)$ $\mid$ $(A,\prec,\succ,\Delta)$ is a dendriform Hopf algebra$\}$.
\item $\LI(A)=\left\{(*,I)\mid \begin{array}{l}
(A,*,\Delta)\mbox{ is a Hopf algebra and}\\
I\mbox{ is a left ideal of }A \mbox{ such that }A_+=Prim(A) \oplus I
\end{array}\right\}$.
\end{enumerate}
There is a bijection between these two sets, given by:
$$\Phi_A: \left\{ \begin{array}{rcl}
\DD(A)&\longrightarrow & \LI(A)\\
(\prec,\succ) &\longrightarrow &(\prec+\succ,A_+ \prec A_+).
\end{array}\right.$$
\end{theo}

\begin{proof} By proposition \ref{30}, $\Phi_A$ is well-defined. We now define the inverse bijection $\Psi_A$.
Let $(*,I) \in \LI(A)$. In order to lighten the notation, we put $V=Prim(A)$.
By corollary \ref{28}, there exists a isomorphism of coalgebras $\phi_I:T(V)\longrightarrow A$, such that $\phi_I(v)=v$ for all $v\in V$
and $\phi_I(T(V)_{\geq 2})=I$. Let $\tilde{*}$ be the product on $T(V)$, making $\phi_I$ an isomorphism of Hopf algebras.
Then $T(V)_{\geq 2}$ is a left ideal of $T(V)$. By proposition \ref{31}, there exists a dendriform structure $(T(V),\tilde{\prec},\tilde{\succ})$ on $T(V)$.
Let $(A,\prec,\succ)$ be the dendriform Hopf algebra structure on $A$, making $\phi_I$ an isomorphism of dendriform Hopf algebras.
We then put $\Psi_A(*,I)=(\prec,\succ)$. Note that $\prec+\succ=*$, as $\tilde{\prec}+\tilde{\succ}=\tilde{*}$.\\

Let us show that $\Phi_A \circ \Psi_A=Id_{\LI(A)}$. Let $(*,I) \in \LI(A)$. We put $\Phi_A \circ \Psi_A(*,I)=(*',I')$.
Then, as $\Phi_I:(T(V),\tilde{\prec},\tilde{\succ})\longrightarrow (T(V),\prec,\succ)$ is an isomorphism of dendriform algebras:
$$I=\Phi_I(T(V)_{\geq _2})=\phi_I(T(V)^{\tilde{\prec} 2})=A^{\prec 2}.$$
By definition of $\Phi_A$, $I'=A^{\prec 2}=I$. Moreover, $*'=\prec+\succ=*$. \\

Let us show that $\Psi_A \circ \Phi_A=Id_{\DD(A)}$. Let $(\prec,\succ) \in \DD(A)$. We put $\Psi_A \circ \Phi_A(\prec,\succ)=(\prec',\succ')$,
and $\Phi_A(\prec,\succ)=(*,I)$. As $\phi_I:(T(V),\tilde{\prec},\tilde{\succ})\longrightarrow (A,\prec',\succ')$ is an isomorphism
of dendriform algebras, we have to prove that $\phi_I:(T(V),\tilde{\prec},\tilde{\succ})\longrightarrow (A,\prec,\succ)$ is an isomorphism of dendriform algebras.
Let $a=u_1\ldots u_k$, $b=v_1\ldots v_l \in T_+(V)$. First:
$$\phi_I(a\tilde{\prec} b)-\phi_I(a)\prec \phi_I(b) \in \phi_I(T_{\geq 2}(V))+A_+\prec A_+=I+I=I.$$
Let us prove that $\phi_I(a\tilde{\prec} b)=\phi_I(a) \prec \phi_I(b)$ by induction on $k$. For $k=1$, let us proceed by induction on $l$. 
For $l=1$, $\phi_I(a\tilde{\prec} b)=\phi_I(v_1u_1)$. Moreover, in $A$:
\begin{eqnarray*}
\tdelta(\phi_I(u_1)\prec \phi_I(v_1))&=&\tdelta(u_1\prec v_1)\\
&=&v_1\otimes u_1,\\
\tdelta(\phi_I(u_1\tilde{\prec} v_1))&=&(\phi_I \otimes \phi_I)\circ \tdelta(u_1\prec v_1)\\
&=&\phi_I(v_1)\otimes \phi_I(u_1)\\
&=&v_1\otimes u_1.
\end{eqnarray*}
So $\phi_I(u_1\tilde{\prec} v_1)-\phi_I(u_1)\prec \phi_I(v_1) \in Prim(A)\cap I=(0)$.  
Let us suppose the result for all $l'<l$. Then, by the induction hypothesis applied to $b''$:
\begin{eqnarray*}
\tdelta(\phi_I(u_1\tilde{\prec} b))&=&\tdelta(\phi_I(bu_1))\\
&=&(\phi_I \otimes \phi_I)(b \otimes u_1+b' \otimes u_1\tilde{\prec} b'')\\
&=&\phi_I(b) \otimes \phi_I(u_1)+\phi_I(b)' \otimes \phi_I(u_1)\tilde \phi_I(b'')\\
&=&\tdelta(\phi_I(u_1)\prec \phi_I(b)).
\end{eqnarray*}
So $\phi_I(u_1\tilde{\prec} b))-\phi_I(u_1) \prec \phi_I(b) \in Prim(A)\cap I=(0)$. This prove the result for $k=1$. Let us assume the result 
at all rank $k'<k$. We put $u_1\ldots u_{k-1}=\tilde{a}$. Then, using the first step:
\begin{eqnarray*}
\phi_I(a\tilde{\prec} b)&=&\phi_I(u_k \tilde{\prec} (\tilde{a}\tilde{*}b))\\
&=&\phi_I(u_k)\prec (\phi_I(\tilde{a})*\phi_I(\tilde{b}))\\
&=&(\phi_I(u_k)\prec \phi_I(\tilde{a}))\prec\phi_I(\tilde{b})\\
&=&(\phi_I(u_k \prec \tilde{a}))\prec\phi_I(\tilde{b})\\
&=&\phi_I(a)\prec\phi_I(\tilde{b}).
\end{eqnarray*}
So $\Psi_A \circ \Phi_A=Id_{\DD(A)}$. \end{proof}\\

%% file: chap4.tex
\section{MuPAD computations}

We here give the different MuPAD procedures we used in this text.

\subsection{Dimensions of simple modules}

The following procedures compute the dimension of the simple module $\Gamma_{(a_1,\ldots,a_n)}$ for $\mathfrak{sl}_n$,
$\mathfrak{sp}_{2n}$, $\mathfrak{so}_{2n}$ and $\mathfrak{so}_{2n+1}$.
\begin{verbatim}
dimsl:=proc(a)
local n,res;
begin
n:=nops(a)+1; 
res:=product(product((sum(a[k],k=i..j-1)+j-i)/(j-i),j=i+1..n),i=1..n-1);
return(res);
end_proc;

dimsp:=proc(a)
local n,k,s,l,res;
begin
n:=nops(a);l:=[];
for k from 1 to n do l:=l.[sum(a[x],x=k..n)+n-k]; end_for; 
res:=product(product((l[i]-l[j])*(l[i]+l[j]+2),j=i+1..n),i=1..n-1);
res:=res*product(l[j]+1,j=1..n)/product((2*n-2*j-1)!,j=0..n-1);
return(res);
end_proc;

dimsoodd:=proc(a)
local n,k,l,res;
begin
n:=nops(a); l:=[];
for k from 1 to n do l:=l.[sum(a[x],x=k..n)-a[n]/2+n-k]; end_for; 
res:=product(product((l[i]-l[j])*(l[i]+l[j]+1),j=i+1..n),i=1..n-1);
res:=res*product(2*l[j]+1,j=1..n)/product((2*n-2*j-1)!,j=0..n-1);
return(res);
end_proc;

dimsoeven:=proc(a)
local n,k,l,res;
begin
n:=nops(a); l:=[];
for k from 1 to n-2 do l:=l.[sum(a[x],x=k..n-2)+a[n-1]/2+a[n]/2+n-k];end_for; 
l:=l.[a[n-1]/2+a[n]/2+1,-a[n-1]/2+a[n]/2];
res:=product(product((l[i]-l[j])*(l[i]+l[j]),j=i+1..n),i=1..n-1);
res:=res/product((2*n-2*j)!,j=1..n-1)*2^(n-1);
return(res);
end_proc;
\end{verbatim}

The following procedures give the representations of dimension smaller than $dim(\g)$ for $\g=\mathfrak{sl}_n$,
$\mathfrak{sp}_{2n}$, $\mathfrak{so}_{2n}$ and $\mathfrak{so}_{2n+1}$.

\begin{verbatim}
smallsl:=proc(n)
local ens,prod,i,d,dg;
begin
dg:=n^2-1;
print(Unquoted,"dimension of g: ".expr2text(dg));
ens:=[];
for i from 1 to n-1 do ens:=ens.[{0,1,2}]; end_for;
prod:=combinat::cartesianProduct::list(ens[x]$x=1..n-1);
for i from 1 to 3^(n-1) do
   d:=dimsl(prod[i]); 
   if d<=dg then print(Unquoted,"heighest weight ".expr2text(prod[i]).
                 " of dimension ".expr2text(d));
   end_if;
end_for;
end_proc;

smallsp:=proc(n)
local ens,prod,i,d,dg;
begin
dg:=n*(2*n+1);
print(Unquoted,"dimension of g: ".expr2text(dg));
ens:=[];
for i from 1 to n do ens:=ens.[{0,1,2}]; end_for;
prod:=combinat::cartesianProduct::list(ens[x]$x=1..n);
for i from 1 to 3^n do
   d:=dimsp(prod[i]); 
   if d<=dg then print(Unquoted,"heighest weight ".expr2text(prod[i]).
                 " of dimension ".expr2text(d));
   end_if;
end_for;
end_proc;

smallsoodd:=proc(n)
local ens,prod,i,d,dg;
begin
dg:=n*(2*n+1);
print(Unquoted,"dimension of g: ".expr2text(dg));
ens:=[];
for i from 1 to n do ens:=ens.[{0,1,2}]; end_for;
prod:=combinat::cartesianProduct::list(ens[x]$x=1..n);
for i from 1 to 3^n do
   d:=dimsoodd(prod[i]); 
   if d<=dg then print(Unquoted,"heighest weight ".expr2text(prod[i]).
                 " of dimension ".expr2text(d));
   end_if;
end_for;
end_proc;

smallsoeven:=proc(n)
local ens,prod,i,d,dg;
begin
dg:=n*(2*n-1);
print(Unquoted,"dimension of g: ".expr2text(dg));
ens:=[];
for i from 1 to n do ens:=ens.[{0,1,2}]; end_for;
prod:=combinat::cartesianProduct::list(ens[x]$x=1..n);
for i from 1 to 3^n do
   d:=dimsoeven(prod[i]); 
   if d<=dg then print(Unquoted,"heighest weight ".expr2text(prod[i]).
                 " of dimension ".expr2text(d));
   end_if;
end_for;
end_proc;
\end{verbatim}

\subsection{Computations for $\mathfrak{sl}_6$}

\label{s4.2} The procedure \verb|testsl6| produces the matrix used in section \ref{s2.5}, corresponding to the action of $\mathfrak{sl}_6$ over $\Lambda^3(V)$,
where $V$ is the standard representation of $\mathfrak{sl}_6$.
\begin{verbatim}
basis:=[[1,2,3],[1,2,4],[1,2,5],[1,2,6],[1,3,4],[1,3,5],[1,3,6],[1,4,5],[1,4,6],
       [1,5,6],[2,3,4],[2,3,5],[2,3,6],[2,4,5],[2,4,6],[2,5,6],
       [3,4,5],[3,4,6],[3,5,6],[4,5,6]]:

indexbasis:=proc(B)
local i;
begin
for i from 1 to 20 do if B=basis[i] then return(i); end_if; end_for;
end_proc:

actionbasis:=proc(A,B)
local vec,vec1,liste,liste2,coef,numero,i,j,k,n;
begin
liste:=[]; i:=B[1]; j:=B[2]; k:=B[3];
vec:=matrix(6,1); vec[i]:=1; vec:=A*vec;
for n from 1 to 6 do liste:=liste.[[vec[n],[n,j,k]]]; end_for;
vec:=matrix(6,1); vec[j]:=1; vec:=A*vec;
for n from 1 to 6 do liste:=liste.[[vec[n],[i,n,k]]]; end_for;
vec:=matrix(6,1); vec[k]:=1; vec:=A*vec; 
for n from 1 to 6 do liste:=liste.[[vec[n],[i,j,n]]]; end_for;
liste2:=[];
for n from 1 to nops(liste) do 
   coef:=(liste[n])[1];
   i:=((liste[n])[2])[1]; j:=((liste[n])[2])[2]; k:=((liste[n])[2])[3]; 
   if (coef<>0) then
      if (i<j) and (j<k) then liste2:=liste2.[[coef,[i,j,k]]]; end_if;
      if (i<k) and (k<j) then liste2:=liste2.[[-coef,[i,k,j]]]; end_if;
      if (j<i) and (i<k) then liste2:=liste2.[[-coef,[j,i,k]]]; end_if;
      if (j<k) and (k<i) then liste2:=liste2.[[coef,[j,k,i]]]; end_if;
      if (k<i) and (i<j) then liste2:=liste2.[[coef,[k,i,j]]]; end_if;
      if (k<j) and (j<i) then liste2:=liste2.[[-coef,[k,i,j]]]; end_if;
   end_if;
end_for;
vec:=matrix(20,1);
for n from 1 to nops(liste2) do
   coef:=(liste2[n])[1];
   numero:=indexbasis((liste2[n])[2]);
   vec1:=matrix(20,1); vec1[numero]:=coef; vec:=vec+vec1;  
end_for;
return(vec);
end_proc:

actionvector:=proc(A,vec)
local i,res;
begin
res:=matrix(20,1);
for i from 1 to 20 do res:=res+vec[i]*actionbasis(A,basis[i]); end_for;
return(res);
end_proc:

testsl6:=proc()
local vec,res,i,j,mat;
begin 
res:=[]; vec:=[];
for i from 1 to 20 do vec:=vec.[a[i]]; end_for;
for i from 1 to 5 do
   mat:=matrix(6,6); mat[i,i]:=1; mat[i+1,i+1]:=-1;
   res:=res.actionvector(mat,vec);
end_for;
for i from 1 to 5 do
   for j from i+1 to 6 do
      mat:=matrix(6,6); mat[i,j]:=1;
      res:=res.actionvector(mat,vec);
      mat:=matrix(6,6); mat[j,i]:=1;
      res:=res.actionvector(mat,vec);
   end_for;
end_for;
return(res);
end_proc;
\end{verbatim}

\subsection{Computations for $\mathfrak{g}_2$}

\label{s4.3} The procedure \verb|testg2| produces the matrix used in section \ref{s2.5}, corresponding to the action of $\mathfrak{g}_2$ on $M_{7,2}(\mathbb{C})$.
\begin{verbatim}
H1:=matrix([[1,0,0,0,0,0,0],[0,-1,0,0,0,0,0],[0,0,2,0,0,0,0],[0,0,0,0,0,0,0],
            [0,0,0,0,-2,0,0],[0,0,0,0,0,1,0],[0,0,0,0,0,0,-1]]):
H2:=matrix([[0,0,0,0,0,0,0],[0,1,0,0,0,0,0],[0,0,-1,0,0,0,0],[0,0,0,0,0,0,0],
            [0,0,0,0,1,0,0],[0,0,0,0,0,-1,0],[0,0,0,0,0,0,0]]):
Y1:=matrix([[0,0,0,0,0,0,0],[1,0,0,0,0,0,0],[0,0,0,0,0,0,0],[0,0,1,0,0,0,0],
            [0,0,0,2,0,0,0],[0,0,0,0,0,0,0],[0,0,0,0,0,-1,0]]):
Y2:=matrix([[0,0,0,0,0,0,0],[0,0,0,0,0,0,0],[0,-1,0,0,0,0,0],[0,0,0,0,0,0,0],
            [0,0,0,0,0,0,0],[0,0,0,0,1,0,0],[0,0,0,0,0,0,0]]):
Y3:=-Y1*Y2+Y2*Y1: Y4:=-1/2*(Y1*Y3-Y3*Y1): 
Y5:=1/3*(Y1*Y4-Y4*Y1): Y6:=Y2*Y5-Y5*Y2:
X1:=matrix([[0,1,0,0,0,0,0],[0,0,0,0,0,0,0],[0,0,0,2,0,0,0],[0,0,0,0,1,0,0],
            [0,0,0,0,0,0,0],[0,0,0,0,0,0,-1],[0,0,0,0,0,0,0]]):
X2:=matrix([[0,0,0,0,0,0,0],[0,0,-1,0,0,0,0],[0,0,0,0,0,0,0],[0,0,0,0,0,0,0],
            [0,0,0,0,0,1,0],[0,0,0,0,0,0,0],[0,0,0,0,0,0,0]]):
X3:=X1*X2-X2*X1: X4:=1/2*(X1*X3-X3*X1):
X5:=-1/3*(X1*X4-X4*X1): X6:=-(X2*X5-X5*X2):
g2:=[H1,H2,X1,X2,X3,X4,X5,X6,Y1,Y2,Y3,Y4,Y5,Y6]:

testg2:=proc()
local i,j,res,A,M;
begin
A:=matrix(7,2);
for i from 1 to 7 do for j from 1 to 2 do A[i,j]:=a[i,j]; end_for; end_for;
res:=[];
for i from 1 to 14 do
   M:=g2[i]*A;
   res:=res.linalg::stackMatrix(M[1..7,1],M[1..7,2]); 
end_for;
return(res);
end_proc;
\end{verbatim}

%% file: prelie.bbl
\providecommand{\bysame}{\leavevmode\hbox to3em{\hrulefill}\thinspace}
\providecommand{\MR}{\relax\ifhmode\unskip\space\fi MR }
\providecommand{\MRhref}[2]{%
  \href{http://www.ams.org/mathscinet-getitem?mr=#1}{#2}
}
\providecommand{\href}[2]{#2}
\begin{thebibliography}{10}

\bibitem{Chapoton1}
Fr\'ed\'eric Chapoton, \emph{Alg\`ebres pr\'e-lie et alg\`ebres de {H}opf
  li\'ees \`a la renormalisation}, C. R. Acad. Sci. Paris S\'er. I Math.
  \textbf{332} (2001), no.~8, 681--684.

\bibitem{Chapoton2}
Fr\'ed\'eric Chapoton and Muriel Livernet, \emph{Pre-{L}ie algebras and the
  rooted trees operad}, Internat. Math. Res. Notices \textbf{8} (2001),
  395--408, arXiv:math/0002069.

\bibitem{Connes}
Alain Connes and Dirk Kreimer, \emph{Hopf algebras, {R}enormalization and
  {N}oncommutative geometry}, Comm. Math. Phys \textbf{199} (1998), no.~1,
  203--242, arXiv:hep-th/9808042.

\bibitem{Fulton}
William Fulton and Joe Harris, \emph{Representation theory}, Graduate Texts in
  Mathematics, vol. 129, Springer-Verlag, New York, 1991, A first course,
  Readings in Mathematics.

\bibitem{Gan}
Wee~Liang Gan and Travis Schedler, \emph{The necklace {L}ie coalgebra and
  renormalization algebras}, J. Noncommut. Geom. \textbf{2} (2008), no.~2,
  195--214.

\bibitem{Kreimer}
Dirk Kreimer, \emph{Combinatorics of (pertubative) {Q}uantum {F}ield {T}heory},
  Phys. Rep. \textbf{4--6} (2002), 387--424, arXiv:hep-th/0010059.

\bibitem{Kubo}
F.~Kubo, \emph{Compatible algebra structures of {L}ie algebras}, Ring theory
  2007, World Sci. Publ., Hackensack, NJ, 2009, pp.~235--239.

\bibitem{Loday1}
Jean-Louis Loday, \emph{Dialgebras}, Lecture Notes in Math., vol. 1763,
  Springer, Berlin, 2001.

\bibitem{Loday2}
Jean-Louis Loday and Maria Ronco, \emph{Combinatorial hopf algebras},
  arXiv:0810.0435, 2009.

\bibitem{Loday3}
Jean-Louis Loday and Maria~O. Ronco, \emph{Hopf algebra of the planar binary
  trees}, Adv. Math. \textbf{139} (1998), no.~2, 293--309.

\bibitem{Oudom}
Jean-Michel Oudom and Daniel Guin, \emph{Sur l'alg\`ebre enveloppante d'une
  alg\`ebre pr\'e-{L}ie}, C. R. Math. Acad. Sci. Paris \textbf{340} (2005),
  no.~5, 331--336, arXiv:math/0404457.

\bibitem{Ronco}
Maria Ronco, \emph{A {M}ilnor-{M}oore theorem for dendriform {H}opf algebras},
  C. R. Acad. Sci. Paris S\'er. I Math. \textbf{332} (2001), no.~2, 109--114.

\bibitem{vanderLaan}
Pepijn van~der Laan and Ieke Moerdijk, \emph{Families of {H}opf algebras of
  trees and pre-{L}ie algebras}, Homology, Homotopy Appl. \textbf{8} (2006),
  no.~1, 243--256, arXiv:math/0402022.

\end{thebibliography}
